\documentclass[review,onefignum,onetabnum]{siamart220329}

\usepackage{amsfonts}
\usepackage{graphicx}
\usepackage{epstopdf}
\usepackage{algorithmic}
\usepackage{amssymb} 
\usepackage{tikz}
\usepackage{tkz-euclide}
\usetikzlibrary{patterns}
\ifpdf
  \DeclareGraphicsExtensions{.eps,.pdf,.png,.jpg}
\else
  \DeclareGraphicsExtensions{.eps}
\fi

\usepackage{amsmath}
\usepackage[shortlabels]{enumitem}
\usepackage[symbol]{footmisc}
\newcommand{\ml}{\left(}
\newcommand{\mr}{\right)}
\newcommand{\ul}[1]{\underline{#1}}
\newcommand{\intr}{\text{int}}
\newcommand{\p}{\partial}
\newcommand{\D}{\mathrm{d}}
\newcommand{\ubar}{\bar{u}}
\newcommand{\vbar}{\bar{v}}

\DeclareMathOperator{\spn}{span}

\nolinenumbers
\newcommand\addt[1]{{#1}}
\newcommand\addtm[1]{{#1}}

\newsiamremark{remark}{Remark}
\newsiamremark{hypothesis}{Hypothesis}
\newsiamthm{claim}{Claim}

\headers{Generic Structural Stability for $2 \times 2$ Conservation Laws}{H. K. Tan and A. L. Bertozzi}

\title{Generic Structural Stability for Riemann Solutions to $2 \times 2$ Systems of Hyperbolic Conservation Laws\thanks{Submitted to the editors February 10, 2025.
\funding{This work was supported by the Simons Foundation Math + X award number 510776, ONR grant N00014-23-1-2565, and NSF grant DMS-2407006.} }}

\author{Hong Kiat Tan\thanks{Department of Mathematics, UCLA, Los Angeles, CA 90095, US (\email{maxtanhk@math.ucla.edu}, \email{bertozzi@g.ucla.edu})}
\and 
Andrea L. Bertozzi\footnotemark[2]
}

\usepackage{amsopn}

\makeatletter
\def\cref@override@label@type#1\@nil#2{%
  \expandafter\def\expandafter\cref@currentlabel\expandafter{%
    \expandafter{\cref@currentlabel}{#2}{}}%
}
\makeatother

\begin{document}

\maketitle

\begin{abstract}
This paper presents a proof of generic structural stability for Riemann solutions to $2 \times 2$ system of hyperbolic conservation laws in one spatial variable, without diffusive terms. This means that for almost every left and right state, shocks and rarefaction solutions of the same type are preserved via perturbations of the flux functions, \addt{the left state, and the right state}. The main assumptions for this proof involve standard assumptions on strict hyperbolicity and genuine non-linearity, a technical assumption on directionality of rarefaction curves, and the regular manifold (submersion) assumption motivated by concepts in differential topology. We show that the structural stability of the Riemann solutions is related to the transversality of the Hugoniot loci and rarefaction curves in the state space. The regular manifold assumption is required to invoke a variant of a theorem from differential topology, Thom's parametric transversality theorem, to show the genericity of transversality of these curves. This in turn implies the genericity of structural stability. We then  apply this theorem to two examples: the p-system and a $2 \times 2$ system governing the evolution of gravity-driven monodisperse particle-laden thin films. In particular, we illustrate how one can verify all the above assumptions for the former, and apply the theorem to different numerical and physical aspects of the system governing the latter. 
\end{abstract}




\section{Introduction}

This paper considers the Riemann problem for $2\times 2$
systems of conservation laws.  This is a classical problem in nonlinear PDE.  We study the problem of the \addt{structural} stability and \addt{the genericity of the structural stability} of 
double-wave solutions under perturbations of the flux functions, and \addt{the left and right states}.  \addt{This is the first paper to address the genericity of structural stability, meaning that it holds almost everywhere.} The problem is motivated by a specific application of particle-laden flow \cite{Murisic2013} in which the derived equations have no closed-form expression for the flux functions, and thus must be solved numerically, requiring approximation.  Much of the classical analytic theory for solutions to the Riemann problem concerns the ``weak shock" case \cite{Lax} in which the left and right state are close.  However, real-world problems often involve strong shocks and in the particle-laden flow problem can even involve singular shocks \cite{Liwangshock}.  Rather than focusing on only the particle-laden flow problem, we consider a general $2\times2$ system and show that ideas from differential topology, specifically transversality, can be used to prove general stability results of double-wave solutions of the Riemann problem under perturbations to the flux \addt{and the given states}.  This paper includes both rigorous analysis for general systems and also two different examples - the well-known p-system \cite{Wendroff_1, Wendroff_2} and the particle-laden flow problem \cite{Murisic2011, Murisic2013}. 

Regarding prior work on flux and left and right state perturbations, it is well known these perturbations lead to small changes in the solution in the $L^1$-sense (see \cite{fronttracking} for scalar conservation laws and \cite{L1stability_system,fronttracking} for systems of hyperbolic conservation laws). \addt{However, this does not guarantee the structural stability of Riemann solutions, thereby motivating this paper.}

The most important overlap with our work is the seminal paper on structurally stable Riemann solutions for $2\times 2$ systems of conservation laws \cite{Schecter}, which in turn was inspired by \cite{Furtado_SS}. Their paper uses a viscous approximation leading to traveling wave connections between constant states. They show that structural stability holds for a wide class of weak entropy solutions to the $2 \times 2$ Riemann problem with respect to perturbations to the left and right states and the flux functions. However, they have to impose certain technical non-degeneracy conditions and transversality assumptions that cannot be verified without knowing the intermediate states beforehand. In addition, their proof is based on the structure of traveling waves to the viscous problem, rather than directly addressing structural stability of the inviscid problem.
In addition, a general framework using the concept of a wave manifold is elucidated in \cite{Global_formalism, Azevedo,Topo_Riemann_2024}.
Other notable works on specific classes of flux functions include \cite{Kong_shock, Kong_rare} for quaslinear systems, \cite{SS_fronttracking} for shock fronts generated by the front-tracking algorithm, and \cite{Sen, WeiShen} for specific systems.
Our work develops a framework that applies to more general Riemann problems (most notably without the viscous approximation).  In order to analyze the problem without the traveling wave ODE framework, we must go back to directly analyzing intersections of the Hugoniot loci from the left and right states.  This framework is used in \cite{KeyfitzE&U} to analyze the intersection of Hugoniot loci and rarefaction curves and hence prove the existence and uniqueness of Riemann solutions under the usual assumptions of a strictly hyperbolic and genuine nonlinear system for a certain class of flux functions satisfying a ``half-plane condition".

\addt{There is also another series of works \cite{Majda1983Stability}, \cite{Majda1983Existence}, and \cite{Majda1984Compressible}, and \cite{MajdaThomann1987CPDE} that dived into the structural stability for system of conservation laws. However, these were either restricted to a specific class of equations, such as the one-dimensional quasi-linear wave equations - which can be written as a $2 \times 2$ system of conservation laws, or system of conservation laws in multiple space dimensions. In the book by Majda \cite{Majda1984Compressible}, he remarked that the structural stability for multidimensional shock fronts is actually easier to prove compared to the problem in one spatial dimension, noting that the shock fronts for polytropic gases in several space variables have \textit{stronger} structural stability than do system of conservation laws in a single spatial variable. Furthermore, ``structural stability" in their papers was only defined with respect to smooth perturbations to the Riemann initial data, which is different from perturbing the left and right states. Moreover, they did not discuss flux perturbations. Nevertheless, their conditions required for structural stability involves having a non-zero determinant of a matrix with some combinations of the right eigenvectors of the Jacobian matrix and the position vector of the intermediate states, which is related to the transversality condition that we employ to prove the genericity of structural stability here.}

Results related to genericity for systems of conservation laws also exist in the literature. \addt{The authors in \cite{Schecter} extended their previous work about the diffusive problem to \cite{Schecter2} addressing different codimensions of Riemann solutions. They show that the structurally stable co-dimension zero solutions have a boundary that are co-dimension one solutions and lack structural stability.} Another relevant work \addtm{is} the well-known Schaeffer Regularity Theorem, which states that for almost every initial data in the Schwartz space, the solution is piecewise smooth with a finite number of shock curves \cite{Schaeffer}, which does not extend to systems \cite{Caravenna}. Another recent work analyzes shocks and singular curves for the Cauchy problem for pressureless gases in 2D with generic smooth initial data \cite{Bressan_pressureless}.  \addt{Nonetheless, works analyzing initial data in the Schwartz space exclude Riemann initial data, the subject of this paper.}

In this paper, we consider a general $2 \times 2$ system of conservation laws:
\begin{equation}\label{general2x2}
\begin{cases}
u_t + \ml F_0(u,v)\mr_x &= 0, \\
v_t + \ml G_0(u,v)\mr_x &= 0, \\
\end{cases}
\end{equation}
with $u = u(x,t), v = v(x,t)$ for $(x,t) \in \mathbb{R} \times [0,\infty)$. Let $U \subseteq \mathbb{R}^2$ be an arbitrary measurable open set, and assume that $F_0,G_0 \in C^2(U)$ be given with the following prescribed Riemann initial data:
\begin{equation}\label{Riemanndata}
(u,v)(x,0) = \begin{cases}
(u_l,v_l) &\text{ for } x< 0, \\
(u_r,v_r) &\text{ for } x > 0, \\
\end{cases}
\end{equation}
with $(u_l,v_l), (u_r,v_r) \in U.$ The goal is to show that for \textbf{generic} $(u_l,v_l), (u_r,v_r) \in U$, for \textbf{any} sufficiently small $C^2$ perturbation to the flux functions $F_0$ and $G_0$ on a compact subset of $U$ and perturbations to the left and right states themselves, unique double-wave entropy solutions of the same type are preserved, exemplifying the concept of \textbf{structural stability}. 
For instance, a double shock solution remains as a double shock solution upon perturbation, albeit with a slight change in the shock amplitudes. This property is crucial in arguing that physical behaviors are preserved upon perturbing the flux functions.
The theorem's statement can be found in Theorem \ref{Approxflux_generic}. 

In contrast to Lax’s seminal work \cite{Lax} where he proved the existence and uniqueness of weak solutions with “weak shocks” and corresponding “weak” rarefactions, using the classical Implicit Function Theorem, this paper addresses “strong” shocks and rarefactions in which the left and right state can be far apart. This naturally leads to differential topology arguments involving transversality. 
\addt{Specifically, there is an {\em intermediate state} corresponding to the intersection between relevant Hugoniot loci/rarefaction curves.
This intersection needs to be transverse for structural stability, using an Implicit Function Theorem on Banach Spaces. We then show that transversality at these intermediate states is generic, by invoking the Foliated Parametric Transversality Theorem \ref{foliated_parametric_transversality}, a slight modification to the classical Thom’s Parametric Transversality Theorem.}

The main assumptions of this theorem on the flux functions $F_0,G_0 \in C^2(U)$ are as follows:
\begin{enumerate}[leftmargin = *, label=(\roman*)]
\item $(F_0,G_0)$ forms a \textbf{strictly hyperbolic} system in $U$, that is, the Jacobian matrix $J(u,v;F_0,G_0)$ given by
\begin{equation}\label{jacobian}
J(u,v;F_0,G_0) = \begin{pmatrix} (F_0)_u & (F_0)_v \\ (G_0)_u & (G_0)_v \end{pmatrix}
\end{equation}
possess two distinct real eigenvalues for each $(u,v) \in U$.
\item $(F_0,G_0)$ forms a \textbf{genuinely non-linear} system in $U$, that is, for each eigenvalue $\lambda_k(u,v;F_0,G_0)$  and the corresponding right eigenvector $\mathbf{r}_k(u,v;F_0,G_0)$ of the Jacobian matrix in \eqref{jacobian}, \addt{with the scaling of the eigenvectors discussed later in the paper}, we have
\begin{equation}\label{genuine_nonlinearity}
\nabla \lambda_k(u,v;F_0,G_0) \cdot \mathbf{r}_k(u,v;F_0,G_0) \neq 0.
\end{equation}
\item $(F_0,G_0)$ satisfies a \textbf{graph condition} in $U$, that is, either $(F_0)_v(u,v) \neq 0$ for all $(u,v) \in U$ or $(G_0)_u(u,v) \neq 0$ for all $(u,v) \in U$. Without loss of generality, by interchanging the roles of $u$ and $v$, we can assume that $(F_0)_v(u,v) \neq 0$ for all $(u,v) \in U$.
\item The Hugoniot loci satisfy the \textbf{regular manifold assumption} as described in Definition \ref{Manifold_assumption}.
\end{enumerate}
Note that (i) and (ii) are standard assumptions, and hold for various physical scenarios. Condition (iii) is a technical assumption that holds for many systems of interest, including the two systems for which we will illustrate the application of our theorem. This condition allows us to prove structural stability for $C^2$ flux functions, and is equivalent to allowing us to describe the wave curves using the standard Cartesian coordinate frame in the state space. In addition to the assumptions above, in order for us to conclude that structural stability is generic, we have to invoke a theorem from differential topology. This requires the Hugoniot loci \addt{(which for $2\times 2$ systems are curves in the plane)} to be a manifold except at the given state, which motivates the \textbf{regular manifold assumption} (iv). 

Furthermore, it is worth noting that, in contrast to \cite{Schecter} and \cite{Schecter2} which utilize a vanishing viscosity argument, we have restricted our attention to entropy solutions satisfying the Lax entropy conditions. However, in this restricted class of solutions, we are able to prove a more general result, namely, that structural stability is \textbf{generic} in the left and right states. In addition, it is related to Theorem \ref{Transversality_generic}, where we prove that transversality at \textbf{any} intermediate state is \addt{generic in the left and right states}.

Another remark is that our theorem requires restricting the flux functions on compact subsets. This is necessary for $C^2(K)$ to be a Banach space and to preserve the uniqueness of Riemann solutions upon perturbation. Furthermore, restricting on a compact set warns the user to be cautious of singular shocks corresponding to the intersection of the Hugoniot loci at infinity, such as in \cite{Liwangshock}, and rarefaction-singular shocks corresponding to the intersection of the Hugoniot locus and the rarefaction curve at infinity, such as in \cite{Liwangrareshock}. These structures could easily deform under perturbation to yield a weak solution rather than a distributional solution.

The results in this paper have many applications in a variety of $2 \times 2$ systems of conservation laws. In Section \ref{sec:p-system}, we consider the p-system representing compressible isentropic flow in Lagrangian coordinates, given by
\begin{equation}\label{p-system}
\left\{ \begin{aligned}
u_t + \ml p(v) \mr_x = 0, \\
v_t - u_x = 0,
\end{aligned} \right.
\end{equation}
with $p \in C^2(0,\infty)$ representing pressure \cite{Wendroff_1,Wendroff_2}, and $u$ and $v$ representing the velocity and specific volume respectively. In Section \ref{sec:PLF}, we consider the system of equations describing gravity-driven monodisperse particle-laden thin films below:
\begin{equation}\label{PLF}
\left\{
\begin{aligned}
\p_t h + \p_x \ml h^3 f\ml \frac{h\phi_0}{h} \mr \mr &= 0, \\
\p_t (h\phi_0) + \p_x \ml  h^3 g\ml \frac{h\phi_0}{h} \mr \mr &= 0.
\end{aligned} \right.
\end{equation}
For \eqref{PLF}, $f$ and $g$ are functions of particle volume fraction $\phi_0 \in [0,\phi_m]$, with $\phi_m$ representing the maximum packing fraction, and $h$ representing the height of the fluid \cite{Murisic2011, Murisic2013}. A common theme among these equations is the presence of functions computed from the solutions of the conservation laws, such as $p$ in \eqref{p-system} and $f$ and $g$ in \eqref{PLF}, which are either obtained experimentally \cite{Wendroff_1, Wendroff_2} or are computationally expensive to evaluate \cite{Murisic2011}.

\addt{Some applications of Theorem \ref{Approxflux_generic} are as follows.  
Fix a ``true'' pair of flux functions $(f,g)$ that satisfies the assumptions of Theorem \ref{Approxflux_generic}. If it satisfies the assumption of Theorem \ref{Approxflux_generic}, then structural stability holds \addtm{for almost every} pair of left $(u_l,v_l)$ and right states $(u_r,v_r)$. Pick such a pair of $(u_l,v_l)$ and $(u_r,v_r)$. Then, we may therefore freely approximate $(f,g)$ by any numerically or physically motivated approximation that is close in the prescribed topology, and the Riemann solution for that pair will retain its structure. These approximations include interpolating between known pressures at specified volumes for the p-system, interpolating $(f,g)$ on precomputed values on a sparse grid and/or fit polynomials or any physically motivated ansatz, and differing flux functions $(\tilde f,\tilde g)$ using different classes of models motivated by slightly different physical arguments. Furthermore, this same theorem implies that the physical structures of the solution to the Riemann problem are the same despite small errors in the experimental measurement of the left and right states. More details on this, together with numerical simulations, are available in Section \ref{sec:PLF}.}

This paper is organized as follows. In Section \ref{sec:math_prelim}, we review key concepts to $2 \times 2$ systems of hyperbolic conservation laws. We then state the main results in Section \ref{sec:main_results}. In particular, the main theorem is stated as Theorem \ref{Approxflux_generic}, which states that for generic left and right states, we have structural stability of double-wave entropy solutions. The proof of the theorem is done via two steps. First, we proceed by proving Theorem \ref{Approxflux} in Section \ref{sec:IFTonBS}, which says that transversality is required for structural stability. Next, we prove Theorem \ref{Transversality_generic} in Section \ref{sec:ModifiedThom}, demonstrating that for generic left and right states, we indeed have transversality at the intermediate states. Technical lemmas required for the proof of Theorem \ref{Approxflux} and Theorem \ref{Transversality_generic} are available in Appendix \ref{Proof_of_C1} and \ref{Proof_of_FPT} respectively. We then proceed to elucidate the applications of Theorem \ref{Approxflux_generic} to the p-system \eqref{p-system} in Section \ref{sec:p-system} and gravity-driven monodisperse particle-laden thin films \eqref{PLF} in Section \ref{sec:PLF}. Finally, we conclude with some remarks and future work in Section \ref{sec:future_work}.

\section{Preliminaries}
\label{sec:math_prelim}
In this section, we will denote the flux functions in \eqref{general2x2} as $(F,G)$. 

As the solution to \eqref{general2x2} with data prescribed in \eqref{Riemanndata} is usually discontinuous, we turn to the notion of a weak solution as defined below.
\begin{definition}
We say that $(u,v) \in L^\infty(\mathbb{R} \times [0,\infty);\mathbb{R}^2)$ is a \textbf{weak solution} to \eqref{general2x2} with prescribed Riemann data \eqref{Riemanndata} if for all $(\varphi,\phi) \in C^\infty_c(\mathbb{R} \times [0,\infty);\mathbb{R}^2)$, we have that 
\begin{equation}\label{weakeq}
\begin{aligned}
&\int_0^\infty \int_{\mathbb{R}} \left\{ \begin{pmatrix} u \\ v \end{pmatrix} \cdot \begin{pmatrix} \varphi \\ \phi \end{pmatrix}_t + \begin{pmatrix} F(u,v) \\ G(u,v) \end{pmatrix} \cdot \begin{pmatrix} \varphi \\ \phi \end{pmatrix}_x \right\} (x,t)\D x \D t + \int_{\mathbb{R}} \left\{ \begin{pmatrix} u \\ v \end{pmatrix} \cdot \begin{pmatrix} \varphi \\ \phi \end{pmatrix} \right\}(x,0) \D x \\
&= 0.
\end{aligned}
\end{equation}
\end{definition}

The Riemann problem for the system \eqref{general2x2} is the identification of a weak solution for $t > 0$  with initial condition \eqref{Riemanndata} specified by the left and right states $(u_l,v_l)$ and $(u_r,v_r)$. According to \cite{Lax}, solutions consist of constant states separated by either rarefaction waves or lines of discontinuity. Furthermore, by translational symmetry, it suffices to consider solutions centered at $x = 0$. 

\subsection{Shocks}\label{shocks}
Along a line of discontinuity, the states must satisfy the  \textbf{Rankine-Hugoniot jump condition} for some shock speed $s$:
\begin{equation}\label{RHJump}
\begin{pmatrix} F(u_r,v_r) - F(u_l,v_l) \\ G(u_r,v_r) - G(u_l,v_l) \end{pmatrix} = s\begin{pmatrix} u_r - u_l \\ v_r - v_l\end{pmatrix}.
\end{equation} 
By eliminating the variable $s$, this is equivalent to the single equation
\begin{equation}\label{RHeq1}
\ml F(u_r,v_r) - F(u_l,v_l) \mr (v_r - v_l) - \ml G(u_r,v_r) - G(u_l,v_l)\mr (u_r -  u_l) = 0.
\end{equation}
Equivalently, for a given state $(u_g,v_g) \in U$, we define the \textbf{Hugoniot objective function} $H: U \times U \times \ml C^2(U)\mr^2 \rightarrow \mathbb{R}$ as
\begin{equation}\label{RH}
H(u,v;u_g,v_g,F,G) := \ml F(u,v) - F(u_g,v_g) \mr \cdot (v - v_g) - \ml G(u,v) - G(u_g,v_g)\mr \cdot (u -  u_g),
\end{equation}
then the set of points $(u,v)$ satisfying $H(u,v;u_g,v_g,F,G) = 0$ is known as the \textbf{Hugoniot locus} originating from $(u_g,v_g)$, denoted by $\mathcal{H}_{(u_g,v_g)}(F,G)$. For notational convenience, we denote $\mathcal{\dot H}_{(u_g,v_g)}(F,G) = \mathcal{H}_{(u_g,v_g)}(F,G) \setminus \{(u_g,v_g)\}$ to be the Hugoniot loci punctured at the given state $(u_g,v_g)$ and  $H_{(u_g,v_g)}(u,v;F,G) := H(u,v;u_g,v_g,F,G)$. The dependence on the flux functions $(F,G)$ will be dropped for brevity, and we will treat $u_g,v_g,F,$ and $G$ as parameters or arguments of the relevant functions when the context is clear. 

For mathematical well-posedness and physical relevance, it is customary to impose an additional entropy condition. Here, we will impose the \textbf{Lax entropy conditions} \cite{Lax}, which states that for arbitrary left and right states given by $(u_l,v_l)$ and $(u_r,v_r)$ respectively along a line of discontinuity, there exists a $k \in \{1,2\}$ such that the propagation speed $s$ in \eqref{RHJump} satisfies:
\begin{equation}\label{LaxEntropy}
\begin{aligned}
\lambda_k(u_l,v_l;F,G) &> s > \lambda_k(u_r,v_r;F,G), \\ 
\lambda_{k-1}(u_l,v_l;F,G) &< s < \lambda_{k+1}(u_r,v_r;F,G). \\
\end{aligned}
\end{equation}
In \eqref{LaxEntropy} above, the inequality is inactive if the subscript is neither $1$ nor $2$. 

Such a line of discontinuity in spacetime is a $k-$\textbf{shock} wave if it satisfies the Rankine-Hugoniot jump condition and the Lax entropy condition for the particular value of $k$ in \eqref{LaxEntropy}. In view of the entropy condition, one can define $\mathcal{H}^-_{k,(u_g,v_g)}(F,G)$ to be the locus in the state space consisting of states $(u,v) \in H^{-1}_{k,(u_g,v_g)}(\cdot;F,G)(0)$ connecting to a given state $(u_g,v_g)$ such that 
$$\begin{aligned}
\lambda_k(u_g,v_g;F,G) > s > \lambda_k(u,v;F,G), \\
\lambda_{k-1}(u_g,v_g;F,G) < s < \lambda_{k+1}(u,v;F,G)
\end{aligned}$$
Similarly, one can define $\mathcal{H}_{k,(u_g,v_g)}^{+}(F,G)$ to be the locus consisting of states $(u,v) \in H^{-1}_{k,(u_g,v_g)}(\cdot;F,G)(0)$ connecting to a given state $(u_g,v_g)$ satisfying
$$\begin{aligned}
\lambda_k(u_g,v_g;F,G) < s < \lambda_k(u,v;F,G), \\
\lambda_{k+1}(u_g,v_g;F,G) > s > \lambda_{k-1}(u,v;F,G).
\end{aligned}$$
The positive/negative sign on the superscript represents connecting states with a higher/lower $k-$eigenvalue than the given state $(u_g,v_g)$. By \eqref{LaxEntropy}, any intermediate states $(u^*,v^*) \notin  \{(u_l,v_l),(u_r,v_r)\}$ for the $2 \times 2$ system must lie on $\mathcal{H}^-_{1,(u_l,v_l)} \cap \mathcal{H}^+_{2,(u_r,v_r)}$.

\subsection{Rarefactions}
Here we consider rarefaction solutions following the derivation in \cite{LeVeque} and \cite{LeFloch} as follows. We look for a non-constant self-similar solution in the similarity variable $\xi := \frac{x}{t}$ of the form $\mathbf{u}(x,t) = \mathbf{w}(\xi)$ for some unknown function $w$ connecting from a given state $(u_g,v_g) \in U$. Plugging this into the conservation laws, we obtain
\begin{equation}\label{rare3}
\ml J(\mathbf{w}(\xi)) - \xi \mr \mathbf{w}'(\xi) = 0,
\end{equation}
with $J(\cdot)$ as the Jacobian matrix evaluated at $(u,v)(\xi) = \mathbf{w}(\xi)$. Hence, we have 
\begin{equation}\label{rare10}
\xi = \lambda_k(\mathbf{w}(\xi))
\end{equation}
and $w(\xi)$ must be parallel to the right eigenvector of the Jacobian matrix corresponding to the eigenvalue $\lambda_k$. In other words, there exists some $\alpha(\xi)$ for each $\xi$ such that $\mathbf{w}'(\xi) = \alpha(\xi)\mathbf{r}_k(w(\xi)).$ This implies that the values of $\mathbf{w}(\xi)$ all lie along an integral curve satisfying the following ODE:
\begin{equation}\label{rare4}
\mathbf{w}'(\xi) = \alpha_k(\xi)\mathbf{r}_k(\mathbf{w}(\xi)), \quad \mathbf{w}(\lambda_k(u_g,v_g)) = (u_g,v_g).
\end{equation}
In texts like \cite{LeVeque} and \cite{LeFloch}, the standard procedure to determine $\alpha_k(\xi)$ involves differentiating    equation \eqref{rare10} to obtain
\begin{equation}\label{rare5}
\alpha_k(\xi) = \frac{1}{\nabla \lambda_k(\mathbf{w}(\xi))\cdot \mathbf{r}_k(\mathbf{w}(\xi))}.
\end{equation}
\addt{Note that generally speaking, the rarefaction solution is single-valued if $\xi$, or equivalently $\lambda_k$, is strictly monotone along the integral curve. In other words, we need $\nabla \lambda_k(\mathbf{w}(\xi))\cdot \mathbf{r}_k(\mathbf{w}(\xi)) \neq 0$ for all $\xi$, which is equivalent to the genuine non-linearity condition in (ii).} 
If the flux functions are only $C^2$, this implies that the right-hand side of \eqref{rare4} would only be continuous in $\xi$, and we cannot guarantee that a unique solution to \eqref{rare4} exists by Picard-Lindel\"of. Instead, one would require the flux functions to be at least $C^3$. On the other hand, we can circumvent this for $C^2$ flux functions by using our graph condition in (iii) as follows. Henceforth, we have the following proposition:
\begin{proposition} ($C^2$-flux well-posedness of graph-rarefactions.)\label{rare6} For each $k \in \{1,2\}$, there exists a unique solution to \eqref{rare4} if the flux functions are $C^2$ and satisfy assumptions (i) to (iv).
\end{proposition}
\textit{Proof of Proposition \ref{rare6}.} We start by deducing that solutions to \eqref{rare4} must exist by Peano existence theorem, for $\xi$ up till when $w(\xi)$ first touches the boundary of $U$. Let $r_k^{(i)}$ be the $i$-th component of the right eigenvector $\mathbf{r}_k$ with eigenvalue $\lambda_k$. For any given $F,G\in C^2(U)$, by assumption (iii), we have 
\begin{equation}\label{rare2}
r_k^{(2)} = -\frac{F_u - \lambda_k}{F_v}r_k^{(1)}.
\end{equation}
This implies that solutions to \eqref{rare4} must satisfy the following ODE:
\begin{equation}\label{rareODE}
\begin{aligned}
\frac{\D}{\D u}v_k(u;u_g,v_g,F,G) &= \Xi_k(u,v_k(u;u_g,v_g,F,G);F,G), \\
v_k(u_g;u_g,v_g,F,G) &= v_g,
\end{aligned}
\end{equation}
with
\begin{equation}\label{rareODE2}
\Xi_k(u,v;F,G) = -\frac{F_u(u,v) - \lambda_k(u,v;F,G)}{F_v(u,v)}.
\end{equation}
Observe that solutions to \eqref{rareODE} is unique up till when the curve touches the boundary of $U$ by Picard-Lindel\"of as the right hand side of \eqref{rareODE} is now $C^1$. The solution to \eqref{rareODE} is known as the $k-$\textbf{rarefaction curve} originating from a given state $(u_g,v_g)$, which we denoted by $\mathcal{R}_{k,(u_g,v_g)}$.  For uniqueness, it suffices to show that the map $\xi \mapsto \mathcal{R}_{k,(u_g,v_g)}$ is bijective. Injectivity follows from the genuine non-linearity condition as $\xi = \lambda_k$ is strictly monotonic along rarefaction curves, while surjectivity follows from the fact that $\mathcal{R}_{k,(u_g,v_g)}$ is precisely the range of $w(\cdot)$. We then conclude that solutions to \eqref{rare4} must be unique.  \hfill $\square$

To complete the description for rarefaction waves, we define a few more terminologies and notations as follows. Along the $k-$rarefaction curves, we denote the portion of the curve, together with the given state $(u_g,v_g)$, for which the corresponding characteristic speed $\lambda_k$ to be increasing and decreasing to be $\mathcal{R}^+_{k,(u_g,v_g)}$ and $\mathcal{R}^-_{k,(u_g,v_g)}$ respectively. By definition, we have $\mathcal{R}_{k,(u_g,v_g)} = \mathcal{R}^+_{k,(u_g,v_g)} \cup \mathcal{R}^-_{k,(u_g,v_g)}$. In addition, by the genuine non-linearity condition, eigenvalues are strictly monotonic along rarefaction curves, and thus $\mathcal{R}^+_{k,(u_g,v_g)} \cap \mathcal{R}^-_{k,(u_g,v_g)} = \{(u_g,v_g)\}$. Last but not least, we can define the $k-$\textbf{rarefaction objective function} $R_k:U \times U \times (C^2(U))^2 \rightarrow \mathbb{R}$ given by
\begin{equation}\label{rareobj}
R_k(u,v;u_g,v_g,F,G) := v - v_k(u;u_g,v_g,F,G).
\end{equation}
For a given $u$, one can view \eqref{rareobj} as the signed difference between the given value of $v$ and that obtained along the integral curve at that value of $u$. Hence, the full rarefaction curves $\mathcal{R}_{k,(u_g,v_g)}$ are the zero level sets of $R_k(\cdot;u_g,v_g,F,G): U \rightarrow \mathbb{R}$.  For notational convenience, we will denote $R_{k,(u_g,v_g)}(u,v) := R_k(u,v;u_g,v_g).$

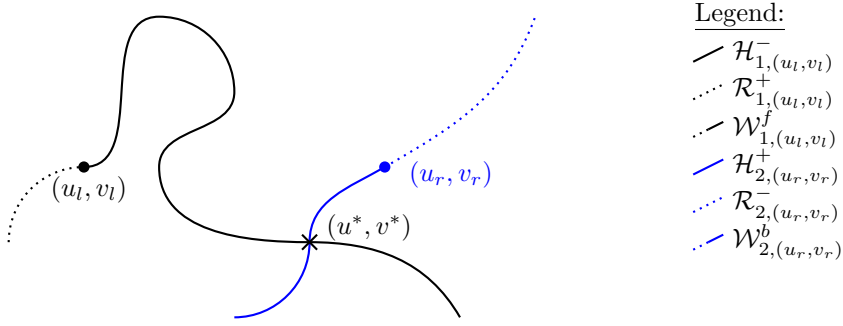
\begin{figure}[ht]
\centering
\begin{tikzpicture}
    \begin{scope}
    \draw[black, thick] (1,1) to[out = 0, in = 180](2,3) to[out = 0, in = 90] (3,2) to[out = -90,in=90] (2,1) to[out = -90, in = 180] (4,0) to [out = 0, in = 120] (6,-1);
    \draw[black, thick, dotted] (0,0) to[out=90,in=180] (1,1);
    \draw[blue, thick] (3,-1) to[out = 0, in = -90] (4,0) to [out = 90, in = -150] (5,1);
    \draw[blue, thick, dotted] (5,1) to[out = 30, in = -110] (7,3);
    \node at (1,1)[circle,fill,inner sep=1.5pt]{};
    \node at (1.7,0.7)[anchor = east]{$(u_l,v_l)$};
    \node at (5,1)[blue, circle,fill,inner sep=1.5pt]{};
    \node at (5.2,0.9)[anchor = west]{  {\color{blue}$(u_r,v_r)$} };
    \draw[black,thick] (3.9,0.1) to (4.1,-0.1);
    \draw[black,thick] (3.9,-0.1) to (4.1,0.1);
    \node at (4.1,0.2)[anchor=west] {$(u^*,v^*)$};
    \node[anchor=west] at (9,3) {\ul{Legend:}};
    \node[anchor=west] at (9.5,2.5) {$\mathcal{H}^-_{1,(u_l,v_l)}$};
    \draw[black, thick] (9.1,2.4) to (9.5,2.6);
    \node[anchor=west] at (9.5,2) {$\mathcal{R}^+_{1,(u_l,v_l)}$};
    \draw[black, thick, dotted] (9.1,1.9) to (9.5,2.1);
    \node[anchor=west] at (9.5,1.5) {$\mathcal{W}^f_{1,(u_l,v_l)}$};
    \draw[black, thick, dotted] (9.1,1.4) to (9.3,1.5);
    \draw[black, thick] (9.3,1.5) to (9.5,1.6);
    \node[anchor=west] at (9.5,1) {$\mathcal{H}^+_{2,(u_r,v_r)}$};
    \draw[blue, thick] (9.1,0.9) to (9.5,1.1);
    \node[anchor=west] at (9.5,0.5) {$\mathcal{R}^-_{2,(u_r,v_r)}$};
    \draw[blue, thick, dotted] (9.1,0.4) to (9.5,0.6);
    \node[anchor=west] at (9.5,0) {$\mathcal{W}^b_{2,(u_r,v_r)}$};
    \draw[blue, thick, dotted] (9.1,-0.1) to (9.3,0);
    \draw[blue, thick] (9.3,0) to (9.5,0.1);
    \end{scope}
\end{tikzpicture}
\caption{A diagram illustrating the geometry of the wave curves and their constituent Hugoniot loci and rarefaction curves with a prescribed left and right state $(u_l,v_l)$ and $(u_r,v_r)$ respectively. We see that since $(u^*,v^*)$ lies along the intersection of the two Hugoniot loci, the solution to the Riemann shock is a double shock problem, with constant states going from $(u_l,v_l)$ to $(u^*,v^*)$ to $(u_r,v_r)$.}
\label{wavecurves1}
\end{figure}

\begin{figure}[ht]
\centering
\begin{tikzpicture}
    \begin{scope}
    \draw[black, thick] (1,1) to[out = 0, in = 180](2,3) to[out = 0, in = 90] (3,2) to[out = -90,in=90] (2,1) to[out = -90, in = 180] (4,0) to [out = 0, in = 120] (6,-1);
    \draw[black, thick, dotted] (0,0) to[out=90,in=180] (1,1);
    \draw[blue, thick] (3,-1) to[out = 0, in = -90] (4,0);
    \draw[blue, thick, dotted] (4,0) to [out = 90, in = -150] (5,1) to[out = 30, in = -110] (7,3);
    \node at (1,1)[circle,fill,inner sep=1.5pt]{};
    \node at (1.7,0.7)[anchor = east]{$(u_l,v_l)$};
    \node at (4,0)[blue, circle,fill,inner sep=1.5pt]{};
    \draw[black,thick] (3.9,0.1) to (4.1,-0.1);
    \draw[black,thick] (3.9,-0.1) to (4.1,0.1);
    \node at (4.1,0.2)[anchor=west] {$(u^*,v^*) = $ {\color{blue}$(u_r,v_r)$}};    
    \node[anchor=west] at (9,3) {\ul{Legend:}};
    \node[anchor=west] at (9.5,2.5) {$\mathcal{H}^-_{1,(u_l,v_l)}$};
    \draw[black, thick] (9.1,2.4) to (9.5,2.6);
    \node[anchor=west] at (9.5,2) {$\mathcal{R}^+_{1,(u_l,v_l)}$};
    \draw[black, thick, dotted] (9.1,1.9) to (9.5,2.1);
    \node[anchor=west] at (9.5,1.5) {$\mathcal{W}^f_{1,(u_l,v_l)}$};
    \draw[black, thick, dotted] (9.1,1.4) to (9.3,1.5);
    \draw[black, thick] (9.3,1.5) to (9.5,1.6);
    \node[anchor=west] at (9.5,1) {$\mathcal{H}^+_{2,(u_r,v_r)}$};
    \draw[blue, thick] (9.1,0.9) to (9.5,1.1);
    \node[anchor=west] at (9.5,0.5) {$\mathcal{R}^-_{2,(u_r,v_r)}$};
    \draw[blue, thick, dotted] (9.1,0.4) to (9.5,0.6);
    \node[anchor=west] at (9.5,0) {$\mathcal{W}^b_{2,(u_r,v_r)}$};
    \draw[blue, thick, dotted] (9.1,-0.1) to (9.3,0);
    \draw[blue, thick] (9.3,0) to (9.5,0.1);
    \end{scope}
\end{tikzpicture}
\caption{A diagram illustrating a single shock solution. Since $(u_r,v_r)$ lies along the Hugoniot locus $\mathcal{H}^-_{1,(u_l,v_l)}$, then $(u^*,v^*) = (u_r,v_r)$. In this case, we will have a single shock solution, with constant states going from $(u_l,v_l)$ to $(u_r,v_r)$.}
\label{wavecurves2}
\end{figure}
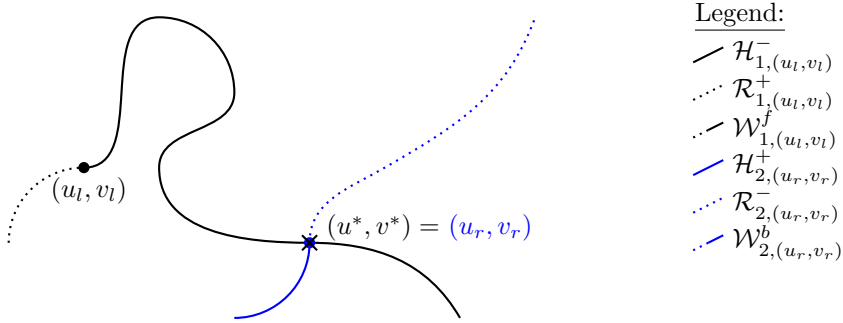

We are now ready to define the notion of an entropy solution as follows:
\begin{definition}
A weak solution to the system of conservation laws \eqref{general2x2} is an \textbf{entropy solution} if all discontinuities occurring in the solution are shocks.
\end{definition}
\addt{Note that the shocks in the above definition refer to the $k-$shocks as described in Section \ref{shocks}.} Next, following the notation in \cite{Liu_shocks}, we define the notion of a wave curve as follows.
\begin{definition}
The \textbf{forward} $k$-\textbf{wave curve} $\mathcal{W}^f_{k,(u_g,v_g)}(F,G)$ and the \textbf{backward} $k$-\textbf{wave curve} $\mathcal{W}^b_{k,(u_g,v_g)}(F,G)$ for a given state $(u_g,v_g)$ are defined by
\begin{equation}\label{wavecurveforward}
\mathcal{W}^f_{k,(u_g,v_g)}(F,G) := \mathcal{R}^+_{k,(u_g,v_g)}(F,G) \cup \mathcal{H}^-_{k,(u_g,v_g)}(F,G), 
\end{equation}
\begin{equation}\label{wavecurvebackward}
\mathcal{W}^b_{k,(u_g,v_g)}(F,G) := \mathcal{R}^-_{k,(u_g,v_g)}(F,G) \cup \mathcal{H}^+_{k,(u_g,v_g)}(F,G).
\end{equation}
\end{definition}
To construct an entropy solution, we first consider all possible intersections between the forward $1-$wave curve $\mathcal{W}^f_{1,(u_l,v_l)}$ originating from the left state and the backward $2-$wave curve $\mathcal{W}^b_{2,(u_r,v_r)}$ originating from $(u_r,v_r)$. We say that an entropy solution exists if the aforementioned curves intersect, and such a solution is unique if the point of intersection is unique. For these unique entropy solutions, denote the point of intersection by $(u_*,v_*)$ (which could be $(u_l,v_l)$ or $(u_r,v_r)$). The corresponding types of solution obtained depending on the types of intersection between the curves are detailed in Table \ref{SolutionMatrix}. Figures \ref{wavecurves1} and \ref{wavecurves2} illustrate some examples of the scenarios described in Table \ref{SolutionMatrix}.

\begin{table}[htbp]
    \centering
\scalebox{0.8}{
    \begin{tabular}{|c|c|c|}
    \hline
     Scenario  &  Type of Solution & States\\
     \hline 
     \multicolumn{3}{|c|}{Single Wave Solutions} \\
    \hline
     $(u_r,v_r) \in \mathcal{H}^-_{1,(u_l,v_l)}$ & Single Shock & $(u_l,v_l)-(u_r,v_r)$ \\
     $(u_r,v_r) \in \mathcal{R}^+_{1,(u_l,v_l)}$  & Single Rarefaction & $(u_l,v_l) \sim (u_r,v_r)$\\
     $(u_l,v_l) \in \mathcal{H}^+_{2,(u_r,v_r)}$  & Single Shock & $(u_l,v_l) - (u_r,v_r)$ \\
     $(u_l,v_l) \in \mathcal{R}^-_{2,(u_r,v_r)}$  & Single Rarefaction & $(u_l,v_l) \sim (u_r,v_r)$\\ 
    \hline 
    \multicolumn{3}{|c|}{Double Wave Solutions ($(u_*,v_*) \neq (u_l,v_l)$ and $(u_*,v_*) \neq (u_r,v_r)$)} \\
    \hline
    $(u^*,v^*) \in \mathcal{H}^-_{1,(u_l,v_l)}$ \& $(u^*,v^*) \in \mathcal{H}^+_{2,(u_r,v_r)}$  & Double Shock & $(u_l,v_l) - (u^*,v^*) - (u_r,v_r)$ \\
    $(u^*,v^*) \in \mathcal{H}^-_{1,(u_l,v_l)}$ \& $(u^*,v^*) \in \mathcal{R}^-_{2,(u_r,v_r)}$  & Shock-Rarefaction & $(u_l,v_l) - (u^*,v^*) \sim (u_r,v_r)$ \\
    $(u^*,v^*) \in \mathcal{R}^+_{1,(u_l,v_l)}$ \& $(u^*,v^*) \in \mathcal{H}^+_{2,(u_r,v_r)}$  & Rarefaction-Shock & $(u_l,v_l) \sim (u^*,v^*) - (u_r,v_r)$ \\
    $(u^*,v^*) \in \mathcal{R}^+_{1,(u_l,v_l)}$ \& $(u^*,v^*) \in \mathcal{R}^-_{2,(u_r,v_r)}$  & Double Rarefaction & $(u_l,v_l) \sim (u^*,v^*) \sim(u_r,v_r)$ \\
    \hline
    \end{tabular}
}
    \caption{A table describing the relationship between the geometry of curves in the state space and the corresponding solution to the Riemann problem. The explicit dependence on $F$ and $G$ for the Hugoniot loci and rarefaction waves are dropped for brevity. When describing the states above, we use $a - b$ to represent a shock which jumps from state $a$ to state $b$, while $a \sim b$ to represent a rarefaction wave connecting from $a$ to $b$ through a continuum of states $s$ between $a$ and $b$.}
    \label{SolutionMatrix}
\end{table}
In view of the definition for a unique entropy solution, we define the following terminologies:
\begin{definition}\label{waveobj}
We say that the unique entropy solution to \eqref{general2x2} is a \textbf{single wave solution} if the resulting solution is either a single shock or a single rarefaction. Otherwise, it is a \textbf{double wave solution}.
\end{definition}



\section{Main Theorems}
\label{sec:main_results}

In this section, we will present the main theorems and their implications. Let $T_x\mathcal{M}$ denote the tangent space to a manifold $\mathcal{M}$ at $x \in \mathcal{M}$. We start by first recalling the notion of transversality of two manifolds (curves):
\begin{definition}\label{transversality_def}
Let $\mathcal{M}$ and $\mathcal{N}$ be submanifolds of a manifold $\mathcal{Y}$. We say that $\mathcal{M}$ and $\mathcal{N}$ are \textbf{transverse} (or equivalently, $\mathcal{M}$ and $\mathcal{N}$ \textbf{intersect transversally}) if  at every point $x \in \mathcal{M} \cap \mathcal{N}$, we have
\begin{equation}\label{T1}
T_x \mathcal{M} + T_x \mathcal{N} = T_x \mathcal{Y}.
\end{equation}
In this case, we write $\mathcal{M} \pitchfork \mathcal{N}.$
\end{definition}
By definition, note that if the two manifolds do not intersect, the transversality property is vacuously satisfied. Next, we define the notion of transversality for the system \eqref{general2x2} with Riemann initial data \eqref{Riemanndata} below:

\begin{definition}\label{Transverse_Def}
Consider the system \eqref{general2x2} with Riemann initial data \eqref{Riemanndata} satisfying assumptions (i), (ii), (iii), and (iv), and the corresponding full Hugoniot loci and full rarefaction connections from the left and right state given by $\mathcal{H}_{(u_l,v_l)}$, $\mathcal{R}_{1,(u_l,v_l)}$, $\mathcal{H}_{(u_r,v_r)}$, and $\mathcal{R}_{2,(u_r,v_r)}$. 
Let $\mathcal{W}_l = \mathcal{H}_{(u_l,v_l)}$ or $\mathcal{R}_{1,(u_l,v_l)}$ and $\mathcal{W}_r = \mathcal{H}_{(u_r,v_r)}$ or $\mathcal{R}_{2,(u_r,v_r)}$ be the left and right full wave curves. We say that the system with initial left and right states $(u_l,v_l), (u_r,v_r)$ a given pair of flux functions $(F_0,G_0)$ possesses the \textbf{transversality property} if the intersection of $\mathcal{W}_l$ and $\mathcal{W}_r$, exclusive of $(u_l,v_l)$ and $(u_r,v_r)$, consists of isolated points $(\ubar,\vbar)$ satisfying
\begin{equation}\label{Transverse_point}
\mathcal{W}_l \pitchfork \mathcal{W}_r.
\end{equation}
Suppose that the flux functions $(F_0,G_0)$ are restricted to some compact subset $K$ of $U$ containing $(u_l,v_l)$ and $(u_r,v_r)$ in its interior. We say that the system possesses the \textbf{transversality property on $K$} if the intersection of $\mathcal{W}_l$ and $\mathcal{W}_r$ restricted to $K$, exclusive of  $(u_l,v_l)$ and $(u_r,v_r)$,\footnote[3]{The restriction of rarefaction curves on a compact set $K$ refers to solving the rarefaction ODE \eqref{rareODE} up till when it first hits the boundary of the compact set $K$. \addt{Note that there is no need to set such a restriction for the Hugoniot loci as they are just loci of points satisfying an algebraic equation, not an integral curve that is constructed by solving some ODE.}
} results in only finitely many of such discrete points satisfying \eqref{Transverse_point}.
\end{definition}
The geometric interpretation of \eqref{Transverse_point} is that the curves in $K \subset \mathbb{R}^2$ intersect transversely at each intersection point $(\ubar,\vbar)$. \\

First, we have the following theorem: \addt{
\begin{theorem}\label{Approxflux} (Structural Stability for Double-Wave Solutions.) Consider the system \eqref{general2x2} satisfying assumptions (i), (ii), (iii), and (iv), with Riemann initial data \eqref{Riemanndata}, and possessing a unique double-wave entropy solution. Then, let $K \subset U$ be a compact subset of $U$ whose interior contains $(u_l,v_l)$, $(u_r,v_r)$  and $(u^*,v^*)$, and suppose that $(F_0,G_0)$ satisfy the transversality property on $K$. Furthermore, consider $F_\delta, G_\delta \in C^2(K)$ as perturbations to the flux functions and $(\delta u_l, \delta v_l), (\delta u_r,\delta v_r)$ as perturbations to the left and right states given by
\begin{equation}\label{perturbed_f&g}
\begin{aligned}
\tilde{F}(u,v) &= (F_0 + F_\delta)(u,v) \\
\tilde{G}(u,v) &= (G_0 + G_\delta)(u,v) 
\end{aligned} \quad \text{ and } \quad 
\begin{aligned}
(\tilde{u}_l, \tilde{v}_l) &= (u_l + \delta u_l, v_l + \delta v_l) \\
(\tilde{u}_r, \tilde{v}_r) &= (u_r + \delta u_r, v_r + \delta v_r).
\end{aligned}
\end{equation}
Then, there exists $\varepsilon_1,\varepsilon_2 > 0$ such that if $\|(F_\delta,G_\delta,\delta u_l, \delta v_l, \delta u_r, \delta v_r)\|_{C^2(K)^2 \times \mathbb{R}^4} < \varepsilon_1$, then the corresponding perturbed $2 \times 2$ system of conservation laws and initial data
\begin{equation}\label{perturbed_system}
\begin{cases}
u_t + \ml \tilde{F}(u,v)\mr_x &= 0, \\
v_t + \ml \tilde{G}(u,v)\mr_x &= 0, \\
\end{cases} \quad \text{ and } \quad (u,v)(x,0) = \begin{cases}
(\tilde{u}_l,\tilde{v}_l) &\text{ for } x< 0, \\
(\tilde{u}_r,\tilde{v}_r) &\text{ for } x > 0, \\
\end{cases}
\end{equation}
satisfies the transversality property on the same compact set $K$ and admits a unique double-wave entropy solution of the same type, with intermediate state 
$(\tilde{u}^*,\tilde{v}^*) \in \intr(K)$ satisfying
\begin{equation}\label{perturbed_int_states}
\|(\tilde{u}^*,\tilde{v}^*) - (u^*,v^*)\|_2 < \varepsilon_2.
\end{equation}
\end{theorem}}

Next, we define what regular value means below:
\begin{definition}
Let $f: \mathcal{M} \rightarrow \mathcal{N}$ be a $C^r$ map with $r \geq 1$ for $C^r$ (with $r \geq 1$) submanifolds $\mathcal{M} \subset \mathbb{R}^m$ and $\mathcal{N} \subset \mathbb{R}^n$ for some positive integers $m$ and $n$. We say that $c \in \mathcal{M}$ is a \textbf{regular value} of a function $f$ if 
for each $x \in f^{-1}(\{c\})$, $x$ is a \textbf{regular point}, that is,
$df_x : T_x \mathcal{M} \rightarrow T_{f(x)} \mathcal{N}$ is surjective. A point $x \in \mathcal{M}$ is a \textbf{critical point} if it is not a regular point.
\end{definition}
\addt{In the definition above, we employ the definition of submanifolds as in \cite{LeeManifold}.} We then introduce the regular manifold assumption below:
\begin{definition}\label{Manifold_assumption}
We say that the system \eqref{general2x2} satisfies the \textbf{regular manifold assumption} if for any $(u_g,v_g) \in U$, the Hugoniot objective function $H_{(u_g,v_g)}$ restricted to $U \setminus \{(u_g,v_g)\}$ has $0$ as a regular value.
\end{definition}
This assumption is necessary for us to invoke the following theorem:

\begin{theorem}\label{RVT}
(Regular Value Theorem/Pre-image Theorem \cite{Hirsch}.) Let $f: \mathcal{M} \rightarrow \mathcal{N}$ be a $C^r$ map for $r \geq 1$ between $C^r$ manifolds $\mathcal{M}$ and $\mathcal{N}$. If $y \in f(\mathcal{M})$ is a regular value, then $f^{-1}(y)$ is a $C^r$ submanifold of $\mathcal{M}$ with dimension $\dim(\mathcal{M}) - \dim(\mathcal{N})$.
\end{theorem}
In other words, this can be applied to prove that a given object is a (sub)manifold by viewing it as a level set of a sufficiently smooth function. In fact, by the regular manifold assumption above, for a given state $(u_g,v_g)$, the punctured Hugoniot loci 
$\mathcal{\dot H}_{(u_g,v_g)} = H_{(u_g,v_g)}|_{U\setminus \{(u_g,v_g)\}}^{-1}(0) $
are manifolds by the regular value theorem with $F_0, G_0 \in C^2(U)$. The full Hugoniot loci could then be computed by
$\mathcal{H}_{(u_g,v_g)} = \mathcal{\dot H}_{(u_g,v_g)} \cup \{(u_g,v_g)\}$. In the language of differential topology, the regular manifold assumption is equivalent to saying that $H_{(u_g,v_g)}:U \setminus \{(u_g,v_g)\} \rightarrow \mathbb{R}$ is a \textbf{submersion} if we restrict its domain to the Hugoniot locus excluding $(u_g,v_g)$. To illustrate that it is necessary to exclude the given state from the set $U$, observe that the map $\ml dH_{(u_g,v_g)} \mr_{(u,v)}:T_{(u,v)} U \rightarrow T_{H_{(u_g,v_g)}(u,v)} \mathbb{R}$ has the following matrix representation
\begin{equation}\label{Manifold_assumption_1}
\ml dH_{(u_g,v_g)} \mr_{(u,v)} = \begin{pmatrix}
A & B
\end{pmatrix}
\end{equation}
with
\begin{equation}\label{Manifold_assumption_2}
\begin{aligned}
A &= F_u(u,v)\cdot(v - v_g) - G_u(u,v)\cdot(u-u_g) - \ml G(u,v) - G(u_r,v_r)\mr, \\
B &= F_v(u,v)\cdot(v - v_g) - G_v(u,v) \cdot(u-u_g) + \ml F(u,v) - F(u_r,v_r)\mr. \\
\end{aligned}
\end{equation}
Thus, we have
\begin{equation}
\ml dH_{(u_g,v_g)} \mr_{(u_g,v_g)} = \begin{pmatrix}
0 & 0 
\end{pmatrix}
\end{equation}
and fails to be subjective for any given state and flux functions. An example of this would be the closed-figure-eight curve in the Keyfitz-Kranzer system \cite{KeyfitzSingular} with vertex at the given state $(u_g,v_g)$ in which the curve fails to be a manifold. In fact, it is well known that for a $2 \times 2$ system, there will be two branches emanating from the given state $(u_g,v_g)$ corresponding to the $1-$ and $2-$ shock connections, respectively \cite{Liu_shocks}, which in turn implies that the Hugoniot locus cannot ever be a manifold if the given state $(u_g,v_g)$ is included. Nonetheless, our results still applies to the any other systems satisfying the appropriate assumptions and the Keyfitz-Kranzer system as removing the vertex would result in two open loops which are manifolds by themselves.

The following theorem illustrates that the transversality property is generic:
\begin{theorem}\label{Transversality_generic}
(Transversality is Generic.)  \addt{Under the assumptions (i) to (iv)}, for almost every $(u_l,v_l,u_r,v_r) \in U^2$,\footnote[2]{``almost every" here is understood in the sense that the set of $(u_l,v_l,u_r,v_r) \in U^2$ such that the aforementioned property does not hold is of measure zero in $U^2$.} the system \eqref{general2x2} with Riemann initial data \eqref{Riemanndata} satisfying assumptions (i) to (iv)     satisfies the transversality property. In particular, for any compact subset $K \subset U$ containing $(u_l,v_l)$ and $(u_r,v_r)$ in its interior with the restriction of $(F_0,G_0)$ onto $K$, the resulting system satisfies the transversality property on $K$.
\end{theorem}

In the statement of Theorem \ref{Transversality_generic}, recall that a non-intersection also satisfies the transversality property.
As a corollary of Theorems \ref{Approxflux} and \ref{Transversality_generic}, we have the following theorem:
\begin{theorem}\label{Approxflux_generic}
(Generic Structural Stability and Approximation Theorem.) Under the assumptions of Theorem \ref{Approxflux}, for almost every $(u_l,v_l,u_r,v_r) \in U^2$, the conclusion of Theorem \ref{Approxflux} holds.
\end{theorem}

\addtm{Last but not least, the following lemma which will be used in the proofs of Theorems \ref{Approxflux} and  \ref{Transversality_generic} illustrates the connection between the two theorems. In particular, the invertibility of a certain matrix is related equivalent to the transversality of the wave curves.
\begin{lemma}\label{Transverse_matrix}
Let $\mathcal{W}_l$ and $\mathcal{W}_r$ be the left and right wave curves intersecting at some $(u^*,v^*) \in U$, and let $W_l, W_r:U \rightarrow \mathbb{R}$ be the corresponding objective functions. Then, $\mathcal{W}_l \pitchfork \mathcal{W}_r$ holds if and only if the matrix
$\begin{pmatrix}
D_u W_l (u^*,v^*) & D_v W_l (u^*,v^*) \\ 
D_u W_r (u^*,v^*) & D_v W_r (u^*,v^*) \\ 
\end{pmatrix}$
is invertible.
\end{lemma}}

\addt{\textit{Proof of Lemma \ref{Transverse_matrix}.}  By utilizing the objective functions, we can compute 
\begin{equation}\label{Transverse_matrix_1}
\begin{aligned}
T_{(u^*,v^*)} \mathcal{W}_{l} &= \ker\{(dW_{l})_{(u^*,v^*)}\} = \spn \{
\begin{pmatrix}
-D_v W_l(u^*,v^*) & D_u W_l(u^*,v^*) \end{pmatrix} \}, \\
T_{(u^*,v^*)} \mathcal{W}_{r} &= \ker\{(dW_{r})_{(u^*,v^*)}\} = \spn \{\begin{pmatrix}
-D_v W_r(u^*,v^*) & D_u W_r(u^*,v^*)
\end{pmatrix}\}.
\end{aligned}
\end{equation}
Hence, condition \eqref{Transverse_point} is equivalent to requiring that the two vectors in \eqref{Transverse_matrix_1} are linearly independent, which in turn is equivalent to the invertibility of the required matrix. \hfill $\square$}

\addt{Before we proceed to prove the relevant theorems, we summarize the implications of the theorems below.
\begin{enumerate}[leftmargin=20pt]
\item Theorem \ref{Approxflux} states that if we have the transversality property, then the solutions to the Riemann problem are structurally stable with regard to perturbations to the flux functions, left state, and the right state.
\item Theorem \ref{Transversality_generic} states that we indeed have the transversality property for almost every left and right states (that is, generic in the left and right states).
\item Henceforth, Theorem \ref{Approxflux_generic} combine these two theorems above, allowing us to conclude that for almost every left and right state, solutions to the Riemann problem are structurally stable with respect to perturbations to the flux functions, left state, and the right state.
\end{enumerate}}

\addt{With the relevant definitions and theorems described this section, we are now ready to prove Theorem \ref{Approxflux} in the section that follows.}

\section{Implicit Function Theorem on Banach Spaces and Proof of Theorem \ref{Approxflux}}\label{sec:IFTonBS}

Before we begin the proof of Theorem \ref{Approxflux}, we first start off by recalling a couple of definitions and theorems on Fr\'echet derivatives as follows. Let $(X,\|\cdot\|_X)$, $(Y,\|\cdot\|_Y)$, and $(Z,\|\cdot\|_Z)$ be Banach spaces with respect to their respective norms. Denote the set of bounded linear maps from $X$ to $Y$ as $\mathcal{B}(X,Y)$. We then have the following definition for Fr\'echet derivative:
\begin{definition}\label{FrechetDerivative}
A map $f: U \subset X \rightarrow Y$ whose domain $U$ is an open subset of a Banach space $X$ and whose range is a Banach space $Y$ is \textit{Fr\'echet differentiable} at $x \in U$ if there is a bounded linear map $A: X \rightarrow Y$ such that
\begin{equation}
\lim_{\|h\|_X \rightarrow 0}\frac{\|f(x+h) - f(x) - Ah\|_Y}{\|h\|_X} = 0.
\end{equation}
\end{definition}
If such a linear map exists, then it is unique, and we write it as \addtm{$A = df(x)$}. If $f$ is Fr\'echet differentiable at each point $x \in U$, then $\addtm{df}: U \rightarrow \mathcal{B}(X,Y)$
is a well-defined map that assigns each point in $x \in U$ to a bounded linear map $\addtm{df}(x):X\rightarrow Y$ that approximates $f$ near $x$. With that, we can define the notion of a continuously differentiable ($C^1$) map between Banach spaces as follows:
\begin{definition}
A map $f: U \subset X \rightarrow Y$ is \textbf{continuously differentiable} ($C^1$) at $x \in U$ if the map $\addtm{df}: U \rightarrow \mathcal{B}(X,Y)$ is continuous at $x$, with $U$ endowed with the norm on $X$ and $\mathcal{B}(X,Y)$ endowed with the operator norm. Furthermore, we say that $f$ is continuously differentiable on $U$ if it is continuously differentiable at each point $x \in U$.
\end{definition}
We now proceed to define partial Fr\'echet derivatives as follows:
\begin{definition}
Suppose that $f: X \times Y \rightarrow Z$ with $X \times Y$ endowed with the norm $\|(x,y)\|_{X \times Y} = \|x\|_X + \|y\|_Y$ for each $(x,y) \in X \times Y$. Fix any $y \in Y$ and consider the map $g:X \rightarrow Z$ given by $g(x) := f(x,y)$. The \textbf{partial Fr\'echet derivative} of $f$ with respect to $x$ is given by $D_xf(x,y) := \addtm{dg}(x).$
\end{definition}
\addtm{Note that in the above definitions, we use the differential notation $df$ to denote the total (Fr\'echet) derivative of a function $f$, while we reserve the notation $D_xf$ for the partial (Fr\'echet) derivative with respect to the argument $x$. Later on in the next section when we write $df_a$ for a point $a$ in the domain, this denotes the differential $df$ evaluated at $a$, that is, $df_a = df(a)$. The argument for which we take the partial derivative should be clear from the context. We adopt this convention to maintain consistency with the usual notations used in differential topology while clearly distinguishing between total and partial derivatives. Furthermore,} it is not difficult to generalize the above concept to a finite product of Banach spaces with an equivalent product topology. Next, we present the following theorem:
\begin{theorem}\label{jointC1}
Let $Y$ be a given Banach space. Consider a finite sequence of Banach spaces $X_1, X_2, \cdots, X_n$, and non-empty open subsets $U_1, U_2, \cdots, U_n$, and let $f:U_1 \times \cdots \times U_n \rightarrow Y$. Denote $D_if$ to be the partial Fr\'echet derivative in the direction of $X_i$. All partial derivatives $D_if: U_1 \times \cdots \times U_n \rightarrow \mathcal{B}\ml X_i, Y\mr$ exist and are continuous if and only if $f$ is continuously differentiable on $U_1 \times \cdots \times U_n$.
\end{theorem}
By Theorem \ref{jointC1} above, to check that a map $f$ is $C^1$, it suffices to check that each partial Fr\'echet derivative is continuous.

Last but not least, the proof of Theorem \ref{Approxflux} relies heavily on the Implicit Function Theorem on Banach spaces, which we will recall from \cite{HunterApplied} as follows:
\begin{theorem}\label{IFTonBS} (Implicit Function Theorem on Banach Spaces.) Suppose that $X, Y,$ and $Z$ are Banach spaces, and $\mathbf{F}:O \subset X \times Y \rightarrow Z$ is a $C^1$ map defined on an open subset $O$ of $X \times Y$. If $(x_0,y_0) \in O$ is a point such that $\mathbf{F}(x_0,y_0) = 0$, and $D_x\mathbf{F}(x_0,y_0):X \rightarrow Z$ is a one-to-one, onto, and bounded linear map, then there exists $\varepsilon_1,\varepsilon_2 > 0$, an open ball of radius $\varepsilon_1$ about $y_0$ denoted by $B_{\varepsilon_1}(y_0) \subset Y$, an open neighborhood $V(x_0)$ of $x_0$, and a unique map $M: B_{\varepsilon_1}(y_0) \rightarrow V(x_0) \subset B_{\varepsilon_2}(x_0)$ such that
\begin{equation}
 \mathbf{F}(M(y),y) = 0 \quad \text{for all } y \in B_{\varepsilon_1}(y_0).
\end{equation}
The function $M$ is continuously differentiable, and
\begin{equation}
dM(x) = -[D_x\mathbf{F}(M(y),y)]^{-1}D_y\mathbf{F}(M(y),y).
\end{equation}
\end{theorem}

Let $\mathcal{O} \subset \ml C^2(K)\mr^2$ be any open neighborhood of $(F_0,G_0)$ such that assumptions (i) to (iv) continue to hold. Before we begin with the proof of the theorem, we first state a technical lemma below.
\addt{\begin{lemma}\label{C1maps}
The maps $H, R_1, R_2: K \times \mathcal{O} \times U \subset K \times \ml C^2(K)\mr^2 \times U \rightarrow \mathbb{R}$ are $C^1$.
\end{lemma}}

\addt{Note that $K \times \mathcal{O} \times U$ is open in $K \times \ml C^2(K)\mr^2 \times U$ by the induced subspace topology, and $\ml C^2(K)\mr^2$ is a Banach space with the following norm:
$$\|(F,G)\|_{\ml C^2(K)\mr^2} = \|F\|_{C^2(K)} + \|G\|_{C^2(K)} \quad \text{ for each }(F,G) \in \ml C^2(K) \mr^2.$$
The proof of Lemma \ref{C1maps} is in Appendix \ref{Proof_of_C1}.

\addtm{Before we begin with the proof of the theorem, we clarify the notation for the objective functions $H$, $R_1$, and $R_2$, which were originally defined in Section \ref{sec:math_prelim}. As stated in Lemma \ref{C1maps}, these maps have domain $K \times \mathcal{O} \times U$. For the Hugoniot objective function, the full argument list is $H(u,v,F,G,u_g,v_g)$ where $(u,v) \in K$ is the state variable, $(F,G) \in \mathcal{O}$ are the flux functions, and $(u_g,v_g) \in U$ is the given state. For convenience, we also write $H_{(u_g,v_g)}(u,v,F,G) := H(u,v,F,G,u_g,v_g)$. In what follows, we may reorder the arguments of $H$ (and similarly for $R_1$, $R_2$) for clarity in the given context. Furthermore, which variables are treated as parameters versus variables will depend on the context. The same conventions apply to $R_1$ and $R_2$.}

We now begin the proof of Theorem \ref{Approxflux} as follows.} \\

\begin{figure}[ht]
\centering
\scalebox{0.7}{\begin{tikzpicture}
\begin{scope}
    \draw[->,thick] (0,0) -- (12,0)node[anchor = west]{$C^2(K)^2\times \mathbb{R}^4$};
    \draw[thick, dashed] (4.5,-0.2) -- (4.5,3.5);
    \node[] at (4.5,-0.5) {$(F_0,G_0,u_l,v_l,u_r,v_r))$};
    \draw[dashed,thick] (2,1) -- (6,5) -- (6,11) -- (2,7) -- (2,1);
    \draw[rotate around={-20:(4,6)},thick] (4,6) ellipse (1.4cm and 2cm);
    \node[] at (5,8) {$K$};
    \node[] at (4,9.5) {$U$};
    \node at (3,6)[circle,fill,inner sep=1.5pt]{};
    \node[] at (3,5.7) {$(u_l,v_l)$};
    \node at (5,7)[circle,fill,inner sep=1.5pt,blue]{};
    \node[blue] at (5,7.3) {$(u_r,v_r)$};
    \draw[black,thick] (2,7) to[out=-30,in=110] (3,6) to[out = -80,in=180] (4,5.5) to[out=0,in=135] (6,5);
    \draw[blue,thick] (6,11) to[out=-135,in=10] (5,7) to[out=-170,in=70] (4,5.5) to[out=-110,in=0] (2,2);
    \draw[red,thick] (3.9,5.6) to (4.1,5.4);
    \draw[red,thick] (3.9,5.4) to (4.1,5.6);
    \node[red] at (4,5.2) {$(u^*,v^*)$};
    \draw[red,dotted,thick] (4,5.5) to (10,5.5);
    \node at (8.9,6)[circle,fill,inner sep=1.5pt]{};
    \node at (11,7.1)[circle,fill,inner sep=1.5pt,blue]{};
    \draw[black,thick,dashed] (8,7) to[out=-30,in=110] (9,6) to[out = -80,in=180] (10,5.5) to[out=0,in=135] (12,5);
    \draw[black,thick] (8.6,6.3) to[out=-30,in=120] (8.9,6) to[out=-60,in=-140] (10,5.7) to[out=50] (11.3,5.5);
    \node at (8.9,6.6)[]{$(\tilde{u}_l,\tilde{v}_l)$};
    \node at (11,7.3)[blue]{$(\tilde{u}_r,\tilde{v}_r)$};
    \draw[blue,thick] (11.4,7) to[out=150,in = 20] (11,7.1) to[out=-160,in=90] (10,5.7) to[out=-90,in=80] (9.4,4.1);
    \draw[blue,thick,dashed] (12,11) to[out=-135,in=10] (11,7) to[out=-170,in=70] (10,5.5) to[out=-110,in=0] (8,2);
    \draw[rotate around={-20:(10,6)},thick] (10,6) ellipse (1.4cm and 2cm);
    \node[] at (11,8) {$K$};
    \draw[red,thick] (9.9,5.8) to (10.1,5.6);
    \draw[red,thick] (9.9,5.6) to (10.1,5.8);
    \draw[red,thick] (4,5.5) to[out=10,in=160] (7,5.3) to[out=-20,in=180] (10,5.7);
    \node[red] at (7.7,4.9) {$(u^*,v^*) = \mathbf{M}^{(1)}(F,G,\tilde{u}_l,\tilde{u}_r,\tilde{v}_l,\tilde{v}_r)$};
    \node[] at (1,0) {$($};
    \node[] at (11,0) {$)$};
    \fill[pattern=north east lines] (1,-0.15) -- (1,0.15) -- (11,0.15) -- (11,-0.15) -- (1,-0.15);
    \node[] at (6.5,0.4) {$B_{\varepsilon_1^{(1)}}(F_0,G_0,u_l,v_l,u_r,v_r)$};
    \end{scope}
\end{tikzpicture}}
\caption{A diagram illustrating the Implicit Function Theorem on Banach space on the product space $K \times C^2(K)^2 \times \mathbb{R}^4$. For each $(F,G,\tilde{u}_l,\tilde{v}_l,\tilde{u}_r,\tilde{v}_r) \in B_{\varepsilon_1^{(1)}}(F_0,G_0,u_l,v_l,u_r,v_r)$, the black curves represents the forward $1-$wave curve emanating from $(u_l,v_l)$, while the blue curves represents the backward $2-$wave curve originating from $(u_r,v_r)$. The red curve represents the $C^1$ map for the intermediate states $(u^*,v^*) = \mathbf{M}^{(1)}(F,G,\tilde{u}_l,\tilde{v}_l,\tilde{u}_r,\tilde{v}_r)$ obtained upon applying the Implicit Function Theorem on Banach space. The dashed curves for a perturbed pair of flux functions $(F,G)$ correspond to the wave curves when it was at $(F_0,G_0)$ with left and right states $(u_l,v_l)$ and $(u_r,v_r)$.
}
\label{IFTfig1}
\end{figure}
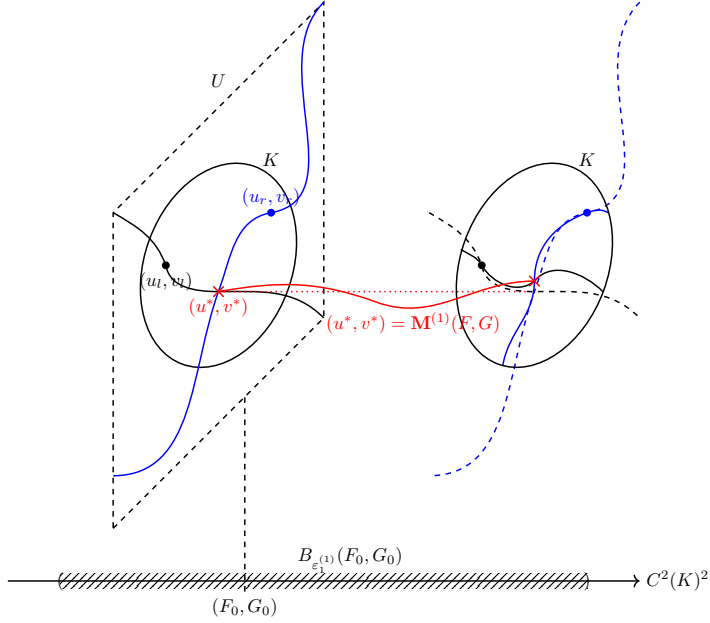

\addt{
\textit{Proof of Theorem \ref{Approxflux}.} The proof can be split into three main steps, mainly, upon perturbation, we have:
\begin{itemize}[leftmargin=30pt]
\item[(A)] Persistence of existence of an intermediate state,
\item[(B)] Persistence of transversality of the corresponding wave curves at that intermediate state, and
\item[(C)] Persistence of uniqueness of such an intermediate state.
\end{itemize}
Without loss of generality, we will assume that the unique intermediate state lies on $\mathcal{H}_{(u_l,v_l)}^+(F_0,G_0)$ and $\mathcal{R}_{2,(u_r,v_r)}^-(F_0,G_0)$. Illustration of the persistence of existence can be found in Figure \ref{IFTfig1}.}\\

\sloppy
\addt{
\ul{(A) Persistence of Existence.} If we define the function $\mathbf{hr}: K \times (C^2(K))^2 \times U^2 \rightarrow \mathbb{R}^2$ given by
\begin{equation}\label{IFT1}
\mathbf{hr}(u,v;F,G,\tilde{u}_l,\tilde{u}_r,\tilde{v}_l,\tilde{v}_r) = \begin{pmatrix} H_{(\tilde{u}_l,\tilde{v}_l)}(u,v,F,G) \\ R_{2,(\tilde{u}_r,\tilde{v}_r)}(u,v,F,G)
\end{pmatrix},
\end{equation}
this implies that $\mathbf{hr}(u^*,v^*;F_0,G_0,u_l,v_l,u_r,v_r) = \mathbf{0}$. Since $\mathbb{R}^2 \times C^2(K)^2 \times \mathbb{R}^4$ is a Banach space, if the hypothesis for Theorem \ref{IFTonBS} holds, we can apply it to $\mathbf{hr}$ to obtain $\varepsilon_1^{(1)}, \varepsilon_2^{(1)} > 0$, an open ball $B_{\varepsilon_1^{(1)}}(F_0,G_0,u_l,v_l,u_r,v_r)$ in $C^2(K)^2 \times \mathbb{R}^4$, an open neighborhood $V^{(1)}(u^*,v^*) \subset B_{\varepsilon_2^{(1)}}(u^*,v^*)$ in $K$, and a unique $C^1$ map $\mathbf{M}^{(1)}: B_{\varepsilon^{(1)}}(F_0,G_0,u_l,v_l,u_r,v_r) \rightarrow V^{(1)}(u^*,v^*)$ with $\mathbf{M}^{(1)}(F_0,G_0,u_l,v_l,u_r,v_r) = (u^*,v^*)$ such that
$$\mathbf{hr}(\mathbf{M}^{(1)}(F,G,\tilde{u}_l,\tilde{v}_l,\tilde{u}_r,\tilde{v}_r),F,G,\tilde{u}_l,\tilde{v}_l,\tilde{u}_r,\tilde{v}_r) = \mathbf{0}$$
for all $(F,G,\tilde{u}_l,\tilde{v}_l,\tilde{u}_r,\tilde{v}_r) \in B_{\varepsilon_1^{(1)}}(F_0,G_0,u_l,v_l,u_r,v_r)$. This implies that the intermediate states for the perturbed pair of flux function $(F,G)$ and states $(\tilde{u}_l,\tilde{v}_l,\tilde{u}_r,\tilde{v}_r)$ given by $(\tilde{u}^*,\tilde{v}^*) = \mathbf{M}^{(1)}(F,G,\tilde{u}_l,\tilde{v}_l,\tilde{u}_r,\tilde{v}_r)$ depends on $(F,G,\tilde{u}_l,\tilde{v}_l,\tilde{u}_r,\tilde{v}_r)$ in a $C^1$-manner. }

\addt{For the solution generated by the Theorem \ref{IFTonBS} to be a valid entropy solution, we have to show that 
\begin{itemize}[leftmargin=20pt]
\item The entropy conditions for the $1$-wave are satisfied by the $1-$shock, and
\item The $2-$eigenvalue increases along the $2-$rarefaction wave.
\end{itemize}
This is because the argument above only shows that $(\tilde{u}^*,\tilde{v}^*) \in  \mathcal{H}_{1,(\tilde{u}_l,\tilde{v}_l)}(F,G) \cap \mathcal{R}_{2,(\tilde{u}_r,\tilde{v}_r)}(F,G)$, but not necessarily in $\mathcal{H}^-_{1,(\tilde{u}_l,\tilde{v}_l)}(F,G) \cap \mathcal{R}^-_{2,(\tilde{u}_r,\tilde{v}_r)}(F,G)$ which constitutes an entropy solution. To do so, let $M^{(1)}_i$ be the $i$-th component of the map $\mathbf{M}^{(1)}$ for $i \in \{1,2\}.$ Consequently, the above properties follow from continuity of the maps $(F,G,\tilde{u}_l,\tilde{v}_l,\tilde{u}_r,\tilde{v}_r) \mapsto \lambda_k(\mathbf{M}^{(1)}(F,G,\tilde{u}_l,\tilde{v}_l,\tilde{u}_r,\tilde{v}_r),F,G,\tilde{u}_l,\tilde{v}_l,\tilde{u}_r,\tilde{v}_r)$, $(F,G,\tilde{u}_l,\tilde{v}_l,\tilde{u}_r,\tilde{v}_r) \mapsto \lambda_k(\tilde{u}_g,\tilde{v}_g,F,G)$, and the shock speed $s: (F,G,\tilde{u}_l,\tilde{v}_l,\tilde{u}_r,\tilde{v}_r) \mapsto \mathbb{R}$ defined by
\begin{equation}\label{IFT2}
s(F,G,\tilde{u}_l,\tilde{v}_l,\tilde{u}_r,\tilde{v}_r) = \begin{cases}
\frac{F(\mathbf{M}^{(1)}(F,G)) - F(\tilde{u}_l,\tilde{v}_l)}{M^{(1)}_1(F,G) - \tilde{u}_l} & \text{ for } u^*(F_0,G_0) \neq u_l, \\
\frac{G(\mathbf{M}^{(1)}(F,G)) - G(\tilde{u}_l,\tilde{v}_l)}{M^{(1)}_2(F,G) - \tilde{v}_l} & \text{ for } v^*(F_0,G_0) \neq v_l, \\
\end{cases}
\end{equation}
from $C^2(K)^2 \times \mathbb{R}^4$ to $\mathbb{R}$ at $(F_0,G_0,u_l,v_l,u_r,v_r)$ for $(\tilde{u}_g,\tilde{v}_g) \in \{(\tilde{u}_l,\tilde{v}_l), (\tilde{u}_r,\tilde{v}_r)\}$ and $k \in \{1,2\}$. Note that these maps are well defined at $(F_0,G_0)$ by the genuine nonlinearity assumption, and continue to be on a sufficiently small $C^2(K)^2\times \mathbb{R}^4$ neighborhood of $(F_0,G_0,u_l,v_l,u_r,v_r)$.} 

\addt{It remains to verify that the hypothesis of Theorem \ref{IFTonBS} holds. By Lemma \ref{C1maps}, one can show that $\mathbf{hr}:K \times (C^2(K))^2 \times U^2 \rightarrow \mathbb{R}^2$ is a $C^1$ map. To verify that $D_{(u,v)}\mathbf{hr}(u^*,v^*;F_0,G_0,u_l,v_l,u_r,v_r): \mathbb{R}^2 \rightarrow \mathbb{R}^2$ is invertible, we do so by first computing the $2 \times 2$ matrix representation of this linear operator as follows:
\begin{equation}\label{IFT3}
D_{(u,v)}\mathbf{hr}(u^*,v^*;F_0,G_0,u_l,v_l,u_r,v_r) = \begin{pmatrix}
D_u H_{(u_l,v_l)} & 
D_v H_{(u_l,v_l)}  \\
D_u R_{2,(u_r,v_r)} & 
D_v R_{2,(u_r,v_r)}  \\
\end{pmatrix}(u^*,v^*,F_0,G_0).
\end{equation}
By the transversality property assumption for the unperturbed pair of flux functions $(F_0,G_0)$ on $K$ and Lemma \ref{Transverse_matrix}, the invertibility of the matrix in \eqref{IFT3} follows. A similar argument holds for a double shock, rarefaction-shock, and a double rarefaction solution, with the corresponding functions defined as
\begin{equation}\label{IFT4}
\begin{aligned}
\mathbf{hh}(u,v;F,G,u_l,v_l,u_r,v_r)  &= \begin{pmatrix} H_{(u_l,v_l)}(u,v,F,G) \\ H_{(u_r,v_r)}(u,v,F,G)
\end{pmatrix},\\
\mathbf{rh}(u,v;F,G,u_l,v_l,u_r,v_r)  &= \begin{pmatrix} R_{1,(u_l,v_l)}(u,v,F,G) \\ H_{(u_r,v_r)}(u,v,F,G)
\end{pmatrix}, \text{ and }\\
\mathbf{rr}(u,v;F,G,u_l,v_l,u_r,v_r)  &= \begin{pmatrix} R_{1,(u_l,v_l)}(u,v,F,G) \\ R_{2,(u_r,v_r)}(u,v,F,G)
\end{pmatrix},\\
\end{aligned}
\end{equation}
respectively.}
\\

\addt{\ul{(B) Persistence of Transversality.} This follows from the continuity of the map $$\det(D_{(u,v)}\mathbf{hr}(\mathbf{M}(\cdot),\cdot)): \ml C^2(K)\mr^2 \times U^2 \rightarrow \mathbb{R}.$$ which in turn implies that the matrix in \eqref{IFT3} continues to be invertible for a sufficiently small $C^2$ perturbation to the flux functions and perturbation to the left and right states. This can be seen from the fact that it only contains terms up to the second derivative of the flux functions $F$ and $G$, $\mathbf{M}^{(1)}, H, R_2$ are $C^1$ maps, and is a well-defined quantity by the genuine non-linearity assumption on a sufficiently small open neighborhood of $(F_0,G_0)$.} \\

\addt{\ul{(C) Persistence of Uniqueness.} 
The above argument only guarantees that the intermediate state $(u^*,v^*)$ obtained for the given pair of flux functions $(F_0,G_0)$ and left and right states $(u_l,v_l,u_r,v_r)$ persists upon perturbation, but does not guarantee that no new intermediate states will surface. Furthermore, it is worth noting that the new intermediate states could be generated by intersecting different left and right wave curves (corresponding to solutions of a different type). To argue that this would not happen for a sufficiently small perturbation, we proceed as follows.}

\addt{By the transversality assumption on the compact set $K$, there are only at most finitely many intersection points between $\mathcal{W}_l$ and $\mathcal{W}_r$, with $\mathcal{W}_l = \mathcal{H}_{(u_l,v_l)}$ or $\mathcal{R}_{1,(u_l,v_l)}$ and $\mathcal{W}_r = \mathcal{H}_{(u_r,v_r)}$ or $\mathcal{R}_{2,(u_r,v_r)}$. Here, we denote the intermediate state corresponding to the unique entropy solution to the system restricted to $K$ in \eqref{general2x2} as $(u^*,v^*)$. For each of these points $i = 1,\cdots, N$ for some positive integer $N$, denote that intersection points as $(\bar{u},\bar{v})^{(i)}$, with $(\bar{u},\bar{v})^{(1)} := (u^*,v^*)$. By applying Theorem \ref{IFTonBS} to each of these intersection points for $i \geq 2$, we generate a finite sequence of $\varepsilon_1^{(i)}, \varepsilon_2^{(i)} > 0$, a finite sequence of open balls $B_{\varepsilon^{(i)}_1}(F_0,G_0,u_l,v_l,u_r,v_r)$, a finite sequence of open neighborhoods $V^{(i)}(\ubar,\vbar)^{(i)}$, and a finite sequence of unique $C^1$ maps $\mathbf{M}^{(i)}:B_{\varepsilon^{(i)}_1}(F_0,G_0,u_l,v_l,u_r,v_r) \rightarrow V^{(i)}(\ubar,\vbar)^{(i)} \subset B_{\varepsilon_{2}^{(i)}}(\ubar,\vbar)^{(i)}$ with $\mathbf{M}^{(i)}(F_0,G_0,u_l,v_l,u_r,v_r) = (\ubar,\vbar)^{(i)}$ such that
$$\mathbf{oo}^{(i)}(\mathbf{M}^{(i)}(F,G,\tilde{u}_l,\tilde{v}_l,\tilde{u}_r,\tilde{v}_r);F,G,\tilde{u}_l,\tilde{v}_l,\tilde{u}_r,\tilde{v}_r) = \mathbf{0}$$
for the appropriate choice of objective functions $\mathbf{oo}^{(i)} \in \{\mathbf{hh},\mathbf{hr},\mathbf{rh},\mathbf{rr}\}$ corresponding to the curves in which $(\ubar,\vbar)^{(i)}$ lie on. By refining the values of $\varepsilon_1^{(i)}$ to be sufficiently small for each $i \geq 2$, we can ensure that these intersection points $(\ubar,\vbar)^{(i)} = \mathbf{M}^{(i)}(F,G,\tilde{u}_l,\tilde{v}_l,\tilde{u}_r,\tilde{v}_r)$ continue to not qualify as entropy solutions for \\$(F,G,\tilde{u}_l,\tilde{v}_l,\tilde{u}_r,\tilde{v}_r) \in B_{\varepsilon_1^{(i)}}(F_0,G_0,u_l,v_l,u_r,v_r)$. Set $\tilde{\varepsilon}_1 = \min_{i \in \{1,\cdots,N\}} \varepsilon_1^{(i)}$ and restrict the domain of each of these maps $\mathbf{M}^{(i)}$ such that we have \\
$\mathbf{M}^{(i)}: B_{\tilde{\varepsilon_1}}(F_0,G_0,u_l,v_l,u_r,v_r) \rightarrow V^{(i)}(\ubar,\vbar)^{(i)}$. These maps remain as the unique maps from $B_{\tilde{\varepsilon}_1}(F_0,G_0,u_l,v_l,u_r,v_r)$ to their corresponding open neighborhoods $V^{(i)}(\ubar,\vbar)^{(i)}$.}

\addt{Last but not least, we need to guarantee that for points in the state space outside of these open neighborhoods $(u,v)$, no new solutions to $\mathbf{oo}^{(i)}(u,v;F,G,\tilde{u}_l,\tilde{u}_r,\tilde{v}_l,\tilde{v}_r) = \mathbf{0}$ for any $\mathbf{oo}^{(i)} \in \{\mathbf{hh},\mathbf{hr},\mathbf{rh},\mathbf{rr}\}$ are generated for a sufficiently small perturbation. Without loss of generality, let us consider all $j \in \{1,\cdots,N\}$ such that $\mathbf{oo}^{(j)} = \mathbf{hh}$. Next, consider the following sets
$$V := \bigcup_{\substack{j = 1, \\ \mathbf{oo}^{(j)} = \mathbf{hh}}} V^{(j)}(\ubar,\vbar)^{(j)}, \quad \text{and} \quad \bar{K} := K \setminus V.$$
Hence, we have that $\mathbf{hh}(u,v;F_0,G_0,u_l,v_l,u_r,v_r) \neq \mathbf{0}$ for each $(u,v) \in \bar{K}$ (else these points would be included at the start of the argument). Since $\bar{K}$ is compact, we have that there exists some constant $M > 0$ such that $\|\mathbf{hh}(u,v;F_0,G_0,u_l,v_l,u_r,v_r)\|_{l^2} \geq M$ uniformly in $(u,v) \in \bar{K}$. 
From Theorem \ref{jointC1} and continuity of $D_{(F,G,\tilde{u}_l,\tilde{u}_r,\tilde{v}_l,\tilde{v}_r)}\mathbf{hh}(u,v;F_0,G_0,u_l,v_l,u_r,v_r): \ml C^2(K)\mr^2 \times \mathbb{R}^4 \rightarrow \mathbb{R}^2$ from Lemma \ref{C1maps}, the Fr\'echet derivative is a bounded operator from $\ml C^2(K) \mr^2 \times \mathbb{R}^4$ to $\mathbb{R}^2$ for $(u,v) \in \bar{K}$. It then follows from the Taylor's Theorem on Banach Spaces that there exists some $\tilde{\varepsilon}_2 > 0$ such that for all $(F,G,\tilde{u}_l,\tilde{u}_r,\tilde{v}_l,\tilde{v}_r) \in B_{\tilde{\varepsilon}_2}(F_0,G_0,u_l,v_l,u_r,v_r)$, we have $\|\mathbf{hh}(u,v;F_0,G_0,u_l,v_l,u_r,v_r)\|_{l^2} \geq \frac{M}{2}$ and hence $\mathbf{hh}(u,v;F_0,G_0,u_l,v_l,u_r,v_r) \neq \mathbf{0}$ for each $(u,v) \in \bar{K}$. By picking $\tilde{\varepsilon} := \min \{ \tilde{\varepsilon}_1,\tilde{\varepsilon}_2 \} > 0$ and restricting the domain of the maps $\mathbf{M}^{(i)}$ to $B_{\tilde{\varepsilon}}(F_0,G_0,u_l,v_l,u_r,v_r)$, no new intermediate states will surface for $(F,G,\tilde{u}_l,\tilde{u}_r,\tilde{v}_l,\tilde{v}_r) \in B_{\tilde{\varepsilon}}(F_0,G_0,u_l,v_l,u_r,v_r)$.
This thus concludes the proof of the theorem.}  \hfill $\square$

\fussy

\section{Foliated Parametric Transversality Theorem and Proof of Theorem \ref{Transversality_generic}}\label{sec:ModifiedThom}

To begin, we recap a few definitions below. Generalizing the notion of transversality between two manifolds, we define the notion of transversality for a sufficiently smooth map $f$ to a manifold as follows.
\begin{definition}\label{transversality_def_2}
Let $f: \mathcal{X} \rightarrow \mathcal{Y}$ be a $C^r$ map with $r \geq 1$, $\mathcal{X}$ and $\mathcal{Y}$ as $C^r$ manifolds, and $\mathcal{Z}$ be a $C^r$ submanifold of $\mathcal{Y}$. We say that $f$ is \textbf{transverse} to $\mathcal{Z}$ if for every $a \in f^{-1}(\mathcal{Z})$, we have
\begin{equation}\label{T3}
df_a(T_a \mathcal{X}) + T_{f(a)}\mathcal{Z} = T_{f(a)}\mathcal{Y}.
\end{equation}
We denote this by $f \pitchfork \mathcal{Z}$.
\end{definition}

\addtm{Recall from Section \ref{sec:IFTonBS} that $df_a$ denotes the differential $df$ evaluated at the point $a$, i.e., $df_a = df(a)$.}

In addition, the notion of transversality between a sufficiently smooth map and a manifold allows us to generalize the regular value theorem in Theorem \ref{RVT} as follows:
\begin{theorem}\label{GPIT}
(Generalized Pre-image Theorem \cite{Hirsch}.) Let $f: \mathcal{X} \rightarrow \mathcal{Y}$ be a $C^r$ map for $r \geq 1$, $\mathcal{X}$ and $\mathcal{Y}$ as $C^r$ manifolds, and $\mathcal{Z}$ be a $C^r$ submanifold of $\mathcal{Y}$. If $f \pitchfork \mathcal{Z}$, then $f^{-1}(\mathcal{Z})$ is a $C^r$ submanifold of $\mathcal{X}$.  Furthermore, the codimension of $f^{-1}(\mathcal{Z})$ in $\mathcal{X}$ is the same as the codimension of $\mathcal{Z}$ in $\mathcal{Y}$.
\end{theorem}

We will now move on to talk about Thom's parametric transversality theorem. This is given as follows:
\begin{theorem}\label{parametric_transversality}
(Thom's Parametric Transversality Theorem \cite{Hirsch}.) Let $\mathcal{X}, \mathcal{P},$ and $\mathcal{Y}$ be $C^r$ manifolds (without boundary) and $\mathcal{Z}$ be a $C^r$ submanifold of \addt{$\mathcal{Y}$}. Consider the maps $F: \mathcal{X} \times \mathcal{P} \rightarrow \mathcal{Y}$ and the associated map $F_p: \mathcal{X} \rightarrow \mathcal{Y}$ for each parameter $p \in \mathcal{P}$.
Suppose that
\begin{enumerate}[leftmargin=*]
\item $r > \max\{0, \dim \mathcal{Y} + \dim \mathcal{Z} - \dim \mathcal{X}\}$,
\item The map $(x,p) \mapsto F_p(x)$ is $C^r$, and
\item $F \pitchfork \mathcal{Z}$.
\end{enumerate}
Then, for almost every $p \in \mathcal{P}$,\footnote{Here, ``almost every" is understood in the sense that the set of $p$ such that this does not hold is of measure $0$ in $\mathcal{P}$.} $F_p \pitchfork \mathcal{Z}$.
\end{theorem}
As described in Section \ref{sec:main_results}, the Hugoniot objective function for a given state $(u_g,v_g) \in U$ has a critical point at $(u,v) = (u_g,v_g)$. Hence, the Hugoniot loci has to be punctured at the given state $(u_g,v_g)$ if we view that as the parameter. This implies that a direct application of Theorem \ref{parametric_transversality} would not be possible, as we require the manifold $\mathcal{X}$ (taking the place of the Hugoniot loci) to be the same for each parameter. This then requires a slight improvement of the parametric transversality theorem as follows:

\begin{theorem}\label{foliated_parametric_transversality} (Foliated Parametric Transversality Theorem.) Let $\mathcal{P}$ and $\mathcal{Y}$ be $C^r$ manifolds, and $\mathcal{Z}$ be a $C^r$ submanifold of $\mathcal{Y}$. Suppose that for each $p \in \mathcal{P}$, we consider a collection of $C^r$ manifolds given by $\{\mathcal{X}_p\}_{p \in \mathcal{P}}$ each with the same dimension $\dim \mathcal{X}$, and the following foliated set:
\begin{equation}\label{XP}
\mathcal{XP} := \bigcup_{p \in \mathcal{P}} \mathcal{X}_p \times \{p\}.
\end{equation}
Consider the maps $F: \mathcal{XP} \rightarrow \mathcal{Y}$ and the associated map $F_p: \mathcal{X}_p \rightarrow \mathcal{Y}$ for each parameter $p \in P$.
Suppose that
\begin{enumerate}[leftmargin=*]
\item $r > \max\{0, \dim \mathcal{Y} + \dim \mathcal{Z} - \dim \mathcal{X}\}$,
\item $\mathcal{XP}$ is a $C^r$ manifold with dimension $\dim \mathcal{XP} = \dim \mathcal{X} + \dim \mathcal{P}$,
\item $T_{(x,p)}\mathcal{XP} \cong T_x \mathcal{X}_p \times T_p \mathcal{P}$ for each $(x,p) \in \mathcal{XP}$,
\item The map $(x,p) \mapsto F_p(x)$ is $C^r$, and
\item $F \pitchfork \mathcal{Z}$.
\end{enumerate}
Then, for almost every $p \in \mathcal{P}$, $F_p \pitchfork \mathcal{Z}$.
\end{theorem}

The proof of this theorem mimics that of Theorem \ref{parametric_transversality} while using assumptions 2. and 3. stated in Theorem \ref{foliated_parametric_transversality} to 
replace arguments utilizing the Cartesian product structure of $\mathcal{X} \times \mathcal{P}$ with $\mathcal{XP}$. The relevant details are available in Appendix \ref{Proof_of_FPT}.

We now proceed with the proof of Theorem \ref{Transversality_generic} as follows. \vspace{3pt}

\textit{Proof of Theorem \ref{Transversality_generic}:} We will do so by applying Theorem \ref{foliated_parametric_transversality}.

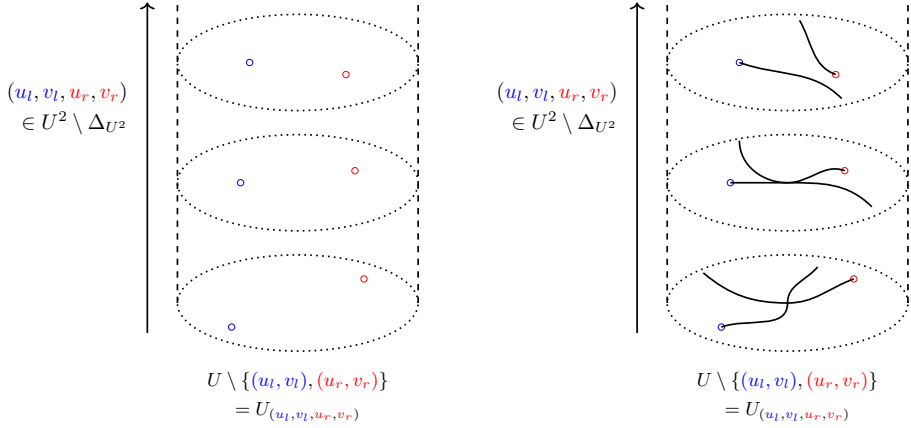
\begin{figure}
\label{Foliated_parametric_fig}
\centering
\scalebox{0.8}{\begin{tikzpicture}
\begin{scope}
\draw[thick,dotted]  (0,0) ellipse (2 and 0.8);
\draw[thick,dotted]  (0,2) ellipse (2 and 0.8);
\draw[thick,dotted]  (0,4) ellipse (2 and 0.8);
\draw[thick,dashed] (-2,0) to (-2,5);
\draw[thick,dashed] (2,0) to (2,5);
\draw[thick,->] (-2.5,-0.5) to (-2.5,5);
\node[anchor=east] at (-2.7,3.5) {$({\color{blue}u_l},{\color{blue}v_l},{\color{red}u_r},{\color{red}v_r}) $};
\node[anchor=east] at (-2.7,3.0) {$\in U^2 \setminus \Delta_{U^2}$};
\node at (0,-1.3) {\small $U\setminus\{{\color{blue}(u_l,v_l)},{\color{red}(u_r,v_r)}\}$};
\node at (0,-1.8) {\small $= U_{({\color{blue}u_l},{\color{blue}v_l},{\color{red}u_r},{\color{red}v_r})}$};
\draw[fill=none,color=red](1.1,0.4) circle (1.5pt);
\draw[fill=none,color=blue](-1.1,-0.4) circle (1.5pt);
\draw[fill=none,color=red](0.95,2.2) circle (1.5pt);
\draw[fill=none,color=blue](-0.95,2.0) circle (1.5pt);
\draw[fill=none,color=red](0.8,3.8) circle (1.5pt);
\draw[fill=none,color=blue](-0.8,4) circle (1.5pt);
\end{scope}
\end{tikzpicture}} \hspace{20pt}
\scalebox{0.8}{\begin{tikzpicture}
\begin{scope}
\draw[thick,dotted]  (0,0) ellipse (2 and 0.8);
\draw[thick,dotted]  (0,2) ellipse (2 and 0.8);
\draw[thick,dotted]  (0,4) ellipse (2 and 0.8);
\draw[thick,dashed] (-2,0) to (-2,5);
\draw[thick,dashed] (2,0) to (2,5);
\draw[thick,->] (-2.5,-0.5) to (-2.5,5);
\node[anchor=east] at (-2.7,3.5) {$({\color{blue}u_l},{\color{blue}v_l},{\color{red}u_r},{\color{red}v_r}) $};
\node[anchor=east] at (-2.7,3.0) {$\in U^2 \setminus \Delta_{U^2}$};
\node at (0,-1.3) {\small $U\setminus\{{\color{blue}(u_l,v_l)},{\color{red}(u_r,v_r)}\}$};
\node at (0,-1.8) {\small $= U_{({\color{blue}u_l},{\color{blue}v_l},{\color{red}u_r},{\color{red}v_r})}$};
\draw[fill=none,color=red](1.1,0.4) circle (1.5pt);
\draw[black,thick](1.1,0.4) to[out=-160,in=0] (0,0) to[out=180,in=-40] (-1.4,0.5);
\draw[fill=none,color=blue](-1.1,-0.4) circle (1.5pt);
\draw[black,thick] (-1.1,-0.4) to[out=20,in=-90] (0,0) to[out=90,in=-130] (0.5,0.6);
\draw[fill=none,color=red](0.95,2.2) circle (1.5pt);
\draw[black,thick](0.95,2.2) to[out=160,in=0] (0,2) to[out=180,in=-90] (-0.8,2.7);
\draw[fill=none,color=blue](-0.95,2.0) circle (1.5pt);
\draw[black,thick] (-0.95,2) to[out=0,in=180] (0,2) to[out=0] (1.4,1.6);
\draw[fill=none,color=red](0.8,3.8) circle (1.5pt);
\draw[black,thick] (0.8,3.8) to[out = 160,in=-60] (0.2,4.7);
\draw[fill=none,color=blue](-0.8,4) circle (1.5pt);
\draw[black,thick] (-0.8,4) to [out = -20] (0.9,3.4);
\end{scope}
\end{tikzpicture}}
\caption{Left: The submanifold $ULR$ of $\mathbb{R}^6$. Right: The punctured Hugoniot loci for each given states ${\color{blue}(u_l,v_l)}$ and ${\color{red}(u_r,v_r)}$ on the doubly-punctured open set $U_{({\color{blue}u_l},{\color{blue}v_l},{\color{red}u_r},{\color{red}v_r})}$ for the double-shock scenario (accounting for Lax entropy conditions). Intuitively, with the $C^1$ dependence of the loci on the given states, we can see that transversality is generic (recall that if the loci do not intersect, they are considered to be transverse too).  }
\end{figure}

\ul{Applying Theorem \ref{foliated_parametric_transversality}.} First, define the punctured set
\begin{equation}\label{U_ulvlurvr}
U_{(u_l,v_l,u_r,v_r)} := U \setminus \{(u_l,v_l),(u_r,v_r)\},
\end{equation}
the diagonal of the set $U^2$ as
\begin{equation}\label{DeltaU^2}
\Delta_{U^2} = \{ (u_l,v_l,u_r,v_r) \in U^2: (u_l,v_l) = (u_r,v_r)\},
\end{equation}
and consider the following set 
\begin{equation}\label{ULR}
ULR := \bigcup_{(u_l,v_l,u_r,v_r) \in U^2 \setminus \Delta_{U^2}} U_{(u_l,v_l,u_r,v_r)} \times \{(u_l,v_l)\} \times \{(u_r,v_r)\}.
\end{equation}
Note that $ULR$ is an open subset of $\mathbb{R}^6$, since it can be written as
\begin{equation}\label{ULR_2}
ULR = U^3 \setminus \ml \{(u,v) = (u_l,v_l)\} \cup \{(u,v) = (u_r,v_r)\} \cup \{(u_l,v_l) = (u_r,v_r)\} \mr.
\end{equation}
Hence, we can endow it with a $C^r$ structure for any $r \geq 1$ so that $ULR$ is a $6-$ dimensional manifold (though not necessarily connected). This is represented on the left of Figure \ref{Foliated_parametric_fig}.

We first apply Theorem \ref{foliated_parametric_transversality} for the double shock scenario. To do so, we consider the function $\mathbf{hh}: ULR \rightarrow \mathbb{R}^2$ given by
\begin{equation}\label{Transverse_H}
\mathbf{hh}(u,v,u_l,v_l,u_r,v_r) = \begin{pmatrix}
H_{(u_l,v_l)}(u,v)\\ H_{(u_r,v_r)}(u,v) \\ 
\end{pmatrix}.
\end{equation}
If we can show that $\mathbf{hh} \pitchfork \{(0,0)\}$, then by the foliated parametric transversality theorem in \eqref{foliated_parametric_transversality}, we deduce that for almost every $(u_l,v_l,u_r,v_r) \in U^2 \setminus \Delta_{U^2}$, the function $\mathbf{hh}_{(u_l,v_l,u_r,v_r)}: U_{(u_l,v_l,u_r,v_r)} \rightarrow \mathbb{R}^2$ defined by
\begin{equation}\label{Transverse_H2}
\mathbf{hh}_{(u_l,v_l,u_r,v_r)}(u,v) = \begin{pmatrix}
H_{(u_l,v_l)}(u,v) \\ H_{(u_r,v_r)}(u,v) 
\end{pmatrix}
\end{equation}
is transverse to $\{(0,0)\}$, ie, $\mathbf{hh}_{(u_l,v_l,u_r,v_r)} \pitchfork \{(0,0)\}.$ Indeed, the assumptions are verified as follows:
\begin{enumerate}[leftmargin=*]
\item $r > \max\{0, \dim \mathbb{R}^2 + \dim \{(0,0)\} - \dim U_{(u_l,v_l,u_r,v_r)}\} = \max\{0, 2 + 0 - 2\} = 0$. Hence, it suffices to pick $r = 1$.
\item $ULR$ is a $C^1$ manifold as explained above, with $\dim ULR = 6 = \dim U_{(u_l,v_l,u_r,v_r)} + \dim (U^2 \setminus \Delta_{U^2})$.
\item $T_{(u,v,u_l,v_l,u_r,v_r)}ULR \cong T_{(u,v)} U_{(u_l,v_l,u_r,v_r)} \times T_{(u_l,v_l,u_r,v_r)} U^2\setminus \Delta_{U^2}$ is true as $ULR$ is just an open subset of $\mathbb{R}^6$ and a $6$-dimensional submanifold of $U$ with a finite number of closed sets of lower dimension removed.
\end{enumerate}

To verify 4. and 5., we compute 
$$d\mathbf{hh}_{(u,v,u_l,v_l,u_r,v_r)}: T_{(u,v,u_l,v_l,u_r,v_r)} ULR \rightarrow T_{\mathbf{hh}(u,v,u_l,v_l,u_r,v_r)} \mathbb{R}^2$$
given by
\begin{equation}\label{T10}
d \mathbf{hh}_{(u,v,u_l,v_l,u_r,v_r)} = \begin{pmatrix} D_u H_{(u_l,v_l)}(u,v) & D_v H_{(u_l,v_l)}(u,v) & A & B & 0 & 0 \\
D_u H_{(u_r,v_r)}(u,v) & D_v H_{(u_r,v_r)}(u,v) & 0 & 0 & C & D\\
\end{pmatrix},
\end{equation}
with
\begin{equation}
\begin{aligned}
A &:=  D_{u_l}H(u,v;u_l,v_l), \quad \\ 
C &:=  D_{u_r}H(u,v;u_r,v_r), \quad
\end{aligned}
\begin{aligned}
B &:=  D_{v_l}H(u,v;u_l,v_l),\\ 
D &:=  D_{v_r}H(u,v;u_r,v_r).
\end{aligned}
\end{equation}
Since each of the entries in the Jacobian computed in \eqref{T10} depends only up to the first derivatives of the flux functions, then $(u,v,u_l,v_l,u_r,v_r) \mapsto \mathbf{hh}_{(u_l,v_l,u_r,v_r)}(u,v)$ can be verified to be $C^1$ map, thus satisfying assumption 4. Next, to verify \eqref{T3} of $\mathbf{hh} \pitchfork \{(0,0)\}$, since $\dim(T_{\mathbf{hh}(u,v,u_l,v_l,u_r,v_r)}\{(0,0)\}) = 0$ and $\dim(T_{\mathbf{hh}(u,v,u_l,v_l,u_r,v_r)}\mathbb{R}^2) = 2$, it suffices to show that $d\mathbf{hh}_{(u,v,u_l,v_l,u_r,v_r)}: T_{(u,v,u_l,v_l,u_r,v_r)}ULR \cong \mathbb{R}^6 \rightarrow \mathbb{R}^2$ surjects onto $\mathbb{R}^2$ for each $(u,v,u_l,v_l,u_r,v_r) \in \mathbf{hh}^{-1}((0,0)).$ \\

To do so, observe that for any $(u_l,v_l,u_r,v_r) \in U^2 \setminus \Delta_{U^2}$, the Hugoniot loci are given by $H^{-1}_{(u_l,v_l)}(0)$ and $H^{-1}_{(u_r,v_r)}(0)$ in $U_{(u_l,v_l,u_r,v_r)}$. Next, consider all $(u,v) \in H^{-1}_{(u_l,v_l)}(0) \cap H^{-1}_{(u_r,v_r)}(0)$ corresponding to intersection points between the Hugoniot loci. By our manifold assumption in Definition \ref{Manifold_assumption}, since $0$ is a regular value for both $H_{(u_l,v_l)}$ and $H_{(u_r,v_r)}$, we have that 
\begin{equation}\label{T9}
\begin{aligned}
(dH_{(u_l,v_l)})_{(u,v)} &= \begin{pmatrix}
D_u H(u,v;u_l,v_l) & D_v H(u,v;u_l,v_l)
\end{pmatrix} \\
(dH_{(u_r,v_r)})_{(u,v)} &= \begin{pmatrix}
D_u H(u,v;u_r,v_r) & D_v H(u,v;u_r,v_r)
\end{pmatrix}\\
\end{aligned}
\end{equation}
both surjects onto $\mathbb{R}$. Observe from the definition of the Hugoniot objective function in \eqref{RH} that 
$(u_l,v_l), (u_r,v_r) \in H^{-1}_{(u,v)}(0)$ for $H_{(u,v)}: U\setminus \{(u,v)\} \rightarrow \mathbb{R}$ since $H(u,v;u_g,v_g) = H(u_g,v_g;u,v)$ for any $(u_g,v_g) \in U$. Physically, this implies that a potential shock connection to $(u,v)$ from a given state $(u_g,v_g)$ is the same as a potential shock connection to $(u_g,v_g)$ by considering the Hugoniot loci emanating from $(u,v)$ instead. The symmetric property described above implies that for each $(u,v,u_g,v_g)$ in $U^2 \setminus \Delta_{U^2}$, we have
\begin{equation}\label{T17}
D_{u}H(u_g,v_g;u,v) = D_{u} H(u,v;u_g,v_g) \text{ and } D_{v}H(u_g,v_g;u,v) = D_{v} H(u,v;u_g,v_g).
\end{equation}
Together with \eqref{T9}, these imply that both 
$$\begin{pmatrix}
A & B
\end{pmatrix} \quad \text{and} \quad \begin{pmatrix}
C & D
\end{pmatrix}$$
surjects onto $\mathbb{R}$. Thus, the map 
$d \mathbf{hh}_{(u,v,u_l,v_l,u_r,v_r)}:\mathbb{R}^6 \rightarrow \mathbb{R}^2$
is surjective for each $(u,v,u_l,v_l,u_r,v_r) \in ULR$, implying that assumption 5. is satisfied and we are done.

To complete the proof, we repeat the arguments above with some modifications for the remaining pairs of curves. We start by defining the foliated and parameterized functions as described in Table \ref{transversality-table1}, and illustrate how to deal with the double-rarefaction scenario, corresponding to $\mathcal{W}_l = \mathcal{R}_{1,(u_l,v_l)}$ and $\mathcal{W}_r = \mathcal{R}_{2,(u_r,v_r)}$.

\begin{table}[htbp]
\begin{center}
\begin{tabular}{|c|c|c|c|}
\hline 
$\mathcal{W}_l$ &  $\mathcal{W}_r$ & Foliated Function & Parametrized Function\\
\hline 
$\mathcal{R}_{1,(u_l,v_l)}$ & $\mathcal{H}_{(u_r,v_r)}$ & $\begin{aligned} \mathbf{hr}(u,v;u_l,v_l,u_r,v_r) \\ = \begin{pmatrix}
R_{1,(u_l,v_l)}(u,v) \\ H_{(u_r,v_r)}(u,v)
\end{pmatrix} \end{aligned}$ & $\begin{aligned} \mathbf{hr}_{(u_l,v_l,u_r,v_r)}(u,v) \\ = \mathbf{hr}(u,v;u_l,v_l,u_r,v_r) \end{aligned}$\\
\hline
$\mathcal{H}_{(u_l,v_l)}$ & $\mathcal{R}_{2,(u_r,v_r)}$ & $\begin{aligned} \mathbf{rh}(u,v;u_l,v_l,u_r,v_r) \\ = \begin{pmatrix}
H_{(u_l,v_l)}(u,v) \\ R_{2,(u_r,v_r)}(u,v)
\end{pmatrix} \end{aligned}$ & $\begin{aligned} \mathbf{rh}_{(u_l,v_l,u_r,v_r)}(u,v) \\ = \mathbf{rh}(u,v;u_l,v_l,u_r,v_r) \end{aligned}$\\
\hline
$\mathcal{R}_{1,(u_l,v_l)}$ & $\mathcal{R}_{2,(u_r,v_r)}$ & $\begin{aligned} \mathbf{rr}(u,v;u_l,v_l,u_r,v_r) \\ = \begin{pmatrix}
R_{1,(u_l,v_l)}(u,v) \\ R_{2,(u_r,v_r)}(u,v)
\end{pmatrix} \end{aligned}$ & $\begin{aligned} \mathbf{rr}_{(u_l,v_l,u_r,v_r)}(u,v) \\ = \mathbf{rr}(u,v;u_l,v_l,u_r,v_r) \end{aligned}$\\
\hline
\end{tabular}
\end{center}
\caption{Choice of foliated and parametrized functions for the remaining cases.}
\label{transversality-table1}
\end{table}

The only key difference is in how we verify assumptions 4 and 5. To do so, we compute
$$d\mathbf{rr}_{(u,v,u_l,v_l,u_r,v_r)}: T_{(u,v,u_l,v_l,u_r,v_r)} ULR \rightarrow T_{\mathbf{rr}(u,v,u_l,v_l,u_r,v_r)} \mathbb{R}^2$$
given by
\begin{equation}\label{T18}
d \mathbf{rr}_{(u,v,u_l,v_l,u_r,v_r)} = \begin{pmatrix} D_u R_{1,(u_l,v_l)}(u,v) & D_v R_{1,(u_l,v_l)}(u,v) & E & F & 0 & 0 \\
D_u R_{1,(u_r,v_r)}(u,v) & D_v R_{2,(u_r,v_r)}(u,v) & 0 & 0 & G & H\\
\end{pmatrix},
\end{equation}
with
\begin{equation}
\begin{aligned}
E &:=  D_{u_l}R_1(u,v;u_l,v_l), \quad \\ 
G &:=  D_{u_r}R_2(u,v;u_r,v_r), \quad 
\end{aligned}
\begin{aligned}
F &:=  D_{v_l}R_1(u,v;u_l,v_l),\\ 
H &:=  D_{v_r}R_2(u,v;u_r,v_r).
\end{aligned}
\end{equation}
Next, we will show that $F$ and $H$ are non-zero, implying that $d\mathbf{rr}_{(u,v,u_l,v_l,u_r,v_r)}:\mathbb{R}^6 \rightarrow \mathbb{R}^2$ is surjective. To see that, for a given $(u_g,v_g)$ and $k \in \{1,2\}$, we start by converting \eqref{rareODE} into the following integral equation:
\begin{equation}\label{T19}
\begin{aligned}
v_k(v;u_g,v_g) = v_g + \int_{u_g}^u\Xi_k(\xi,v_k(\xi;u_g,v_g)) \D \xi.
\end{aligned}
\end{equation}
Taking the derivative with respect to $v_g$, we obtain
\begin{equation}\label{T20}
\frac{\p v_k}{\p v_g}(u;u_g,v_g) = 1 + \int^u_{u_g} \frac{\p \Xi_k}{\p v}(\xi,v_k(\xi;u_g,v_g)) \frac{\p v_k}{\p v_g}(\xi;u_g,v_g) \D \xi.
\end{equation}
Taking the derivative with respect to $u$, we obtain the following ODE for 
\begin{equation}
\begin{cases}
\begin{aligned}
\frac{\p}{\p u}\ml \frac{\p v_k}{\p v_g}\mr (u;u_g,v_g) &= \frac{\p \Xi_k}{\p v}(u,v_k(u;u_g,v_g)) \frac{\p v_k}{\p v_g}(u;u_g,v_g), \\
\ml \frac{\p v_k}{\p v_g}\mr (u_g;u_g,v_g) 
&= 1. 
\end{aligned}
\end{cases}
\end{equation}
This yields
\begin{equation}\label{T21}
\frac{\p v_k}{\p v_g}(v;u_g,v_g) = \exp\ml \int^u_{u_g} \frac{\p \Xi_k}{\p v}(\xi,v_k(\xi;u_g,v_g)) \D \xi \mr.
\end{equation}
Together with the definition of $R_k$ in \eqref{rareobj}, we have
\begin{equation}\label{T22}
F = D_{v_l} R_1(u,v;u_l,v_l) = - \exp\ml \int^u_{u_g} \frac{\p \Xi_k}{\p v}(\xi,v_k(\xi;u_g,v_g)) \D \xi \mr.
\end{equation}
We can easily show that for a strictly hyperbolic pair of $C^2$ flux functions, $\frac{\p \Xi_k}{\p v}$ is a continuous function and hence $F \neq 0$ for any $(u,v) \in U$. This can also be used to justify interchanging the derivative and integral in \eqref{T20} by the Dominated Convergence Theorem for a given $(u,v,u_g,v_g) \in U^2 \setminus \Delta_{U^2}$. An analogous argument would imply that $H \neq 0$ too and we are done. In addition, by performing a similar computation for $E$ and $G$, the astute reader can verify that the following equalities are true:
\begin{equation}\label{T30}
\begin{pmatrix} E & F \end{pmatrix} = F \begin{pmatrix}  -\Xi_1(u_l,v_l) & 1 \end{pmatrix} \text{ and } \begin{pmatrix} G & H \end{pmatrix} = H \begin{pmatrix} -\Xi_2(u_r,v_r) & 1 \end{pmatrix}.
\end{equation}
Geometrically, it is worth noting that \eqref{T30} implies the orthogonality of the vectors $\begin{pmatrix} E & F \end{pmatrix}$ and $\begin{pmatrix} G & H \end{pmatrix}$ to $\mathbf{r}_1(u_l,v_l)$ and $\mathbf{r}_2(u_r,v_r)$ respectively. The remaining cases follow from a similar argument, utilizing the fact that $\begin{pmatrix}A & B \end{pmatrix}, \begin{pmatrix}C & D \end{pmatrix}, \begin{pmatrix}E & F \end{pmatrix}$ and $\begin{pmatrix}G & H \end{pmatrix}$ are surjective onto $\mathbb{R}$.

\ul{Implications and Proving Theorem \ref{Transversality_generic}.}

Recall that by Theorem \ref{foliated_parametric_transversality}, we have shown that $\mathbf{hh}_{(u_l,v_l,u_r,v_r)}\pitchfork\{(0,0)\}$ for almost every $(u_l,v_l,u_r,v_r) \in U^2 \setminus \Delta_{U^2}$. This implies the following:
\begin{itemize}[leftmargin = 20pt]
\item The Hugoniot loci $H^{-1}_{(u_l,v_l)}(0)$ and $H^{-1}_{(u_r,v_r)}(0)$
intersect at discrete points. By the generalized pre-image theorem (Theorem \ref{GPIT}), $\mathbf{hh}_{(u_l,v_l,u_r,v_r)}^{-1}(0,0)$ is a $0$-dimensional submanifold of $U_{(u_l,v_r)}$. Since $$\mathbf{hh}_{(u_l,v_l,u_r,v_r)}^{-1}(0,0) = H^{-1}_{(u_l,v_l)}(0) \cap H^{-1}_{(u_r,v_r)}(0),$$ upon restricting to a compact set $K$ containing $(u_l,v_l)$ and $(u_r,v_r)$ in its interior, the set of discrete intersection points is a closed subset of the compact set $K$ by continuity of $\mathbf{hh}_{(u_l,v_l,u_r,v_r)}$ and hence must be finite.
\item By Definition \ref{transversality_def_2}, we have for each $(u^*,v^*) \in \mathbf{hh}^{-1}_{(u_l,v_l,u_r,v_r)}(0,0),$
\begin{equation}\label{T6}
\ml d\mathbf{hh}_{(u_l,v_l,u_r,v_r)}\mr_{(u^*,v^*)}(T_{(u^*,v^*)} U) + T_{(0,0)}\{(0,0)\} = T_{(0,0)}\mathbb{R}^2.
\end{equation}
Since $\dim(T_{(0,0)}\{(0,0)\}) = 0$, $\dim(T_{(0,0)}\mathbb{R}^2) = 2$, and $\ml 
d \mathbf{hh}_{(u_l,v_l,u_r,v_r)}\mr_{(u^*,v^*)}: \mathbb{R}^2 \rightarrow \mathbb{R}^2$, this implies that
\begin{equation}\label{T23}
\begin{aligned}
\ml 
d\mathbf{hh}_{(u_l,v_l,u_r,v_r)}\mr_{(u^*,v^*)} =& \begin{pmatrix}
D_u H_{(u_l,v_l)}(u^*,v^*) & D_v H_{(u_l,v_l)}(u^*,v^*) \\
D_u H_{(u_r,v_r)}(u^*,v^*) & D_v H_{(u_r,v_r)}(u^*,v^*) \\
\end{pmatrix} \\
\end{aligned}
\end{equation}
is surjective and thus invertible. By Lemma \ref{Transverse_matrix}, we have $T_{(u^*,v^*)}\mathcal{H}_{(u_l,v_l)} + T_{(u^*,v^*)}\mathcal{H}_{(u_r,v_r)} = \mathbb{R}^2$ and we are done. A similar argument holds for the remaining cases.
\end{itemize}

Note that the above arguments hold for almost every $(u_l,v_l,u_r,v_r) \in U^2 \setminus_{\Delta_{U^2}}$. Since $\Delta_{U^2}$ is a set of measure zero in $U^2$, the above result thus holds for almost every $(u_l,v_l,u_r,v_r) \in U^2$.
\hfill $\square$

\section{Applications to the p-system}\label{sec:p-system}

Recall that the p-system is given as a collection of two conservation laws is given in \eqref{p-system}
with $U = \mathbb{R} \times (0,\infty)$ and $p \in C^2(0,\infty)$. Here, we have $F_0(u,v) = p(v)$ and $G_0(u,v) = -u$. It can be shown that the system is strictly hyperbolic if $p'(v) < 0$ on $(0,\infty)$, and genuinely nonlinear if $p''(v) > 0$ on $(0,\infty)$. The former holds by thermodynamics arguments \cite{Wendroff_1}, while the latter holds in some cases \cite{Wendroff_2, Bethe}. Furthermore, $(G_0)_u(u,v) = -1 \neq 0$ implies that assumption (iii) is satisfied. To check that the regular manifold assumption holds, we first compute the Hugoniot objective function for $(u,v;u_g,v_g) \in U^2$:
$$H(u,v;u_g,v_g) = (u - u_g)^2 + (p(v) - p(v_g))\cdot (v - v_g),$$
which is a $C^2$ map from $U$ to $\mathbb{R}$ for a given $(u_g,v_g) \in U$ if $p \in C^2(0,\infty)$. Then, for a given $(u_l,v_l,u_r,v_r) \in U^2 \setminus \Delta_{U^2}$, we compute the differential of the map $H_{(u_g,v_g)}$ for $(u_g,v_g) = (u_l,v_l)$ or $(u_r,v_r)$ at any $(u,v) \in H_{(u_g,v_g)}^{-1}(0)  \setminus \{(u_g,v_g)\}$ as follows:
\begin{equation}\label{p-system_1}
(dH_{(u_g,v_g)})_{(u,v)} = \begin{pmatrix}
2(u-u_g) & p'(v)\cdot(v-v_g) + p(v) - p(v_g)
\end{pmatrix}.\end{equation}
Observe that this linear map $(dH_{(u_g,v_g)})_{(u,v)}:\mathbb{R}^2 \rightarrow \mathbb{R}$ is clearly surjective if $u \neq u_g.$ If $u = u_g$ but $v \neq v_g$, then $(u,v) \in H^{-1}_{(u_g,v_g)}(0) \setminus \{(u_g,v_g)\}$ implies that
\begin{equation}\label{p-system_2}
(p(v) - p(v_g))\cdot (v - v_g) = 0.
\end{equation}
Since $v \neq v_g$, we must have $p(v) = p(v_g)$, which in turn implies that
\begin{equation}\label{p-system_3}
(dH_{(u_g,v_g)})_{(u_g,v)} = (0 \quad\; p'(v)\cdot(v-u_g) ).
\end{equation}
As $p' < 0$ and $v \neq v_g$, then $p'(v)\cdot(v-v_g) \neq 0$, which implies that $(dH_{(u_g,v_g)})_{(u,v_g)}:\mathbb{R}^2 \rightarrow \mathbb{R}$ is surjective. Hence, the assumptions of Theorem \ref{Approxflux_generic} hold, implying structural stability for almost every left and right state.

Physically, the $p$-system governs compressible isentropic flow in Lagrangian coordinates, with the function $p$ representing pressure. In many applications, the pressure is only provided as a table of values obtained from experiments, and intermediate values are obtained by interpolation. Wendroff analyzed the structure of solutions when $p$ loses its convexity (which implies that the system loses genuine nonlinearity) upon interpolating between points in a table of specific volume and pressure in \cite{Wendroff_1}. In addition to the analysis performed in \cite{Wendroff_1} and \cite{Wendroff_2}, Theorem \ref{Approxflux_generic} implies that for almost every left and right state, if we pick a compact subset $K$ of $U$, a sufficiently good interpolating approximation to the pressure preserves unique double-wave entropy Riemann solutions of the same type. \addt{On top of that, if there are some slight measurement errors to the velocity and specific volume as a Riemann initial data, we can be sure that the same physical structure of the solution persists.}

\section{Applications to Gravity-Driven Particle-Laden Thin Films}\label{sec:PLF}

In this section, we will look at three applications of Theorem \ref{Approxflux_generic} to gravity-driven particle-laden thin films as follows. In particular, we consider a thin-film flow of a uniform mixture of spherical, monodisperse, noncolloidal, negatively buoyant, rigid particles, along with a viscous liquid, down a wide inclined channel that makes an angle $\alpha$ with the horizontal. We employ a diffusive flux model as described in \cite{Murisic2011} and \cite{Murisic2013}. For generalizations to the model, we refer the readers to \cite{luong2025brazilnuteffectbidisperseparticle} and \cite{2024spiral}. Let $x$ denote the coordinates of the mixture along the direction of the incline, $z$ denote the direction normal to the incline, and $t$ denote the time elapsed. For each $(x,t)$, we let $h(x,t) > 0$ denote the height of the free surface of the mixture, $s(x,t) \in [0,1]$ denote the non-dimensionalized height of the free surface, and $\phi_0(x,t) \in [0,\phi_m]$ denote the $z-$averaged particle volume fraction, with $\phi_m := 0.61$ representing the maximum packing fraction. 

By the conservation of suspension volume and the number of particles, we obtain the following $2 \times 2$ system of hyperbolic conservation laws:
\begin{equation}\label{PLF1}
\begin{aligned}
\p_t h + \p_x \ml F(h,h\phi_0) \mr &= 0, \\
\p_t (h\phi_0) + \p_x \ml G(h,h\phi_0) \mr &= 0,
\end{aligned}
\end{equation}
with the flux functions $F$ and $G$ given by
\begin{equation}\label{PLF2}
\begin{aligned}
F(h,h\phi_0) = h^3 f\ml \frac{h\phi_0}{h}\mr = h^3 f(\phi_0), \quad 
&G(h,h\phi_0) = h^3 g\ml \frac{h\phi_0}{h}\mr = h^3 g(\phi_0), \\
f(\phi_0(x,t)) = \int^1_0 u(x,s,t;\phi_0(x,t)) \D s,  \quad &g(\phi_0(x,t)) = \int^1_0 u\phi(x,s,t;\phi_0(x,t)) \D s. 
\end{aligned}
\end{equation}
with $f(\phi_0)$ and $g(\phi_0)$ as functions depending only on $\phi_0$ and $u$ is a scalar function representing the velocity field defined in \eqref{PLF3}. On physical grounds, we will attempt to solve \eqref{PLF1} for $(h,\phi_0) \in (0,\infty) \times (0,\phi_m)$. The corresponding state space is the open set
\begin{equation}\label{PLF5}
U := \{(h,h\phi_0): (h,\phi_0) \in (0,\infty) \times (0,\phi_m)\}
\end{equation}
corresponding to an open triangle in $\mathbb{R}^2$. We then consider the corresponding Riemann initial data
\begin{equation}\label{PLFID}
(h,h\phi_0)(0,x) = \begin{cases}
(h_l,h_l (\phi_0)_l) &\text{ for } x< 0, \\
(h_r,h_r (\phi_0)_r) &\text{ for } x > 0, \\
\end{cases}
\end{equation}
with $h_l,h_r \in \mathbb{R}$ and $(\phi_0)_l, (\phi_0)_r \in (0,\phi_m)$, and hence $(h_l,h_l(\phi_0)_l), (h_r,h_r(\phi_0)_r) \in U$.

In addition, for any fixed $(x,t)$, the velocity field $u(x,s,t;\phi_0(x,t))$ \addt{in \eqref{PLF2}} at each $s \in [0,1]$ can be obtained from the expression
\begin{equation}\label{PLF3}
u(x,s,t;\phi_0(x,t)) = \mu_l \int^s_0 \sigma(x,s',t;\phi_0(x,t)) \ml 1 - \frac{\phi(x,s',t;\phi_0(x,t))}{\phi_m}\mr^2 \D s'
\end{equation}
with the constant $\mu_l$ representing the viscosity of the liquid. Furthermore, for a fixed $(x,t)$ and hence $\phi_0(x,t) \in [0,\phi_m]$, the shear stress $\sigma(s) := \sigma(x,s,t;\phi_0(x,t))$ and the particle volume fraction $\phi(s) := \phi(x,s,t;\phi_0(x,t))$ as a function of $s \in [0,1]$ are obtained by solving the system of nonlinear ODEs below:
\begin{equation}\label{PLF4}
\begin{aligned}
\phi'(s) &=  \frac{(-B_2 + (B_2 + 1)\phi(s) + \rho_s \phi(s)^2)(\phi_m - \phi(s))}{\sigma(s) (\phi_m + (B_1 - 1)\phi(s))}H(\phi(s))H(\phi_m - \phi(s)), \\
\sigma'(s) &= - 1 - \rho_s \phi(s),\\
\sigma(0) &= 1 + \rho_s \phi_0,\\
\sigma(1) &= 0,\\
\end{aligned}
\end{equation}
for some positive physical constants 
$\rho_s, B_1,$ and $B_2$ with $H(\cdot)$ representing the Heaviside function with $H(0) := 0.$ The functions $f$ and $g$ are smooth except at a single critical point $\phi_0 = \phi_c \in [0,\phi_m]$, corresponding to a phase transition from a settled to a ridged regime \cite{Murisic2011}. Depending on the angle of the incline with the horizontal, $\alpha$, the functions $f$ and $g$ demonstrate varying regularity across $\phi_c$. For instance, at $\alpha = 25^\circ$, it can be shown numerically and asymptotically that $f$ and $g$ are $C^2$ across $\phi_c$.

\subsection{Interpolating Flux Functions}\label{PLF_intp}

From the mathematical setup above, one can observe that evolving the system of conservation laws \eqref{PLF1} in time with a numerical scheme would require solving a nonlinear system of ODEs in \eqref{PLF4}, followed by performing three different numerical integrations as observed in \eqref{PLF2} and \eqref{PLF3}. A typical workaround for this is to precompute the solutions to the system of ODEs in \eqref{PLF4} and hence the velocity field $u$ in \eqref{PLF3} and the functions $f$ and $g$ in \eqref{PLF2}. This is done over a sufficiently fine grid in the parameter space $\phi_0 \in [0,\phi_m]$. To obtain the value of $f$ and $g$ for an arbitrary $\phi_0 \in [0,\phi_m]$, an interpolation over the values of $f$ and $g$ precomputed at these grid points would be utilized. For instance, let $N + 1$ be the number of grid points with $(\phi_0)_i$ uniformly distributed on $[0,\phi_m]$ and $0 = (\phi_0)_0 < (\phi_0)_2 < \cdots < (\phi_0)_N = \phi_m$ with grid size $\Delta \phi_0 = \frac{\phi_m}{N}$. Denote $f_i := f((\phi_0)_i)$ and $g_i := g((\phi_0)_i)$ for $i = 1, \cdots, N$ and consider the cubic spline interpolation $f_{\text{int},N}$ and $g_{\text{int},N}$ for $(\phi_0) \in ((\phi_0)_i,(\phi_0)_{i+1})$ with weights $w_i = \frac{\phi_0 - (\phi_0)_i}{\Delta \phi_0}$ and expression given by
\begin{equation}\label{PLF7}
\begin{aligned}
f_{\text{int},N}(\phi_0) &= \alpha_i f_{i-1} + \beta_i f_i + \gamma_i f_{i+1} + \xi_i f_{i+2} \hfill \text{ for } \phi_0 \in [(\phi_0)_i,(\phi_0)_{i+1}), \\
g_{\text{int},N}(\phi_0) &= \alpha_i g_{i-1} + \beta_i g_i + \gamma_i g_{i+1} + \xi_i g_{i+2} \hfill \text{ for } \phi_0 \in [(\phi_0)_i,(\phi_0)_{i+1}),\\
\alpha_i &= \frac{-2w_i+3w_i^2-w_i^3}{6}, \quad \beta_i = \frac{2 - w_i - 2w_i^2 + w_i^3}{2}, \\
\gamma_i &= \frac{2w_i+w_i^2-w_i^3}{2}, \quad
\xi_i = \frac{-w_i + w_i^3}{6}, \\
\end{aligned}
\end{equation}
and the corresponding $2 \times 2$ system of conservation laws:
\begin{equation}\label{PLF8}
\begin{aligned}
\p_t h + \p_x \ml F_{\text{int},N}(h,h\phi_0) \mr &= 0, \\
\p_t (h\phi_0) + \p_x \ml G_{\text{int},N}(h,h\phi_0) \mr &= 0,
\end{aligned}
\end{equation}
with 
\begin{equation}\label{PLF9}
F_{\text{int},N}(h,h\phi_0) = h^3 f_{\text{int},N}\ml \frac{h\phi_0}{h}\mr, \quad 
G_{\text{int},N}(h,h\phi_0) = h^3 g_{\text{int},N}\ml \frac{h\phi_0}{h}\mr. \\
\end{equation}
For consistency, we impose a numerical boundary condition for a four-point cubic spline interpolation described above by introducing two ghost grid cells with values $f_{-1} := f_0 - 3f_1 + 3f_2$ and $f_{N+1} :=  f_{N-2} - 3f_{N-1} + 3f_N$.
By Theorem \ref{Approxflux_generic}, we obtain the following numerical analysis result as a corollary:

\sloppy
\begin{proposition}\label{PLF_prop1}
\addt{Suppose that $f$ and $g$ are $C^2$ across
the phase transition at $\phi_0 = \phi_c$ and the system satisfies assumptions (i) to (iv). 
For almost every $(h_l,h_l(\phi_0)_l,h_r,h_r(\phi_0)_r) \in U^2$ and any compact set $K$ containing $(h_l,h_l(\phi_0)_l)$, $(h_r,h_r(\phi_0)_r),$ and $(h^*,h^*(\phi_0)^*)$, if there is a unique entropy solution to \eqref{PLF1} and \eqref{PLFID}, then for a sufficiently large positive integer $N$, the unique double-wave entropy solutions to \eqref{PLF1} and \eqref{PLF8} are of the same type with their intermediate states sufficiently close.}
\end{proposition}
\fussy
Before we begin with the proof of Proposition \ref{PLF_prop1}, it is worth noting that prior to this paper, there were no known results justifying the structural stability of Riemann solutions to the system of conservation laws obtained via interpolation. In addition, the proof of the convergence of $f_{\text{int},N}$ and $g_{\text{int},N}$ to $f$ and $g$ in $C^2$ requires a little more work as $f$ and $g$ is only $C^2$ across $\phi_c$, though this can be greatly shortened if $f$ and $g$ are $C^3$ across $\phi_c$.\\ 

\textit{Proof.} Let $\varepsilon_1 > 0$ be given. Since $w_i'(\phi_0) = \frac{1}{\Delta \phi_0}$, we have
\begin{equation}\label{PLF10}
\begin{aligned}
f'_{\text{int},N}(\phi) &=  \alpha_i' f_{i-1} + \beta_i' f_i +  \gamma_i' f_{i+1} + \xi_i' f_{i+2}, \\
f''_{\text{int},N}(\phi) &=  \alpha_i'' f_{i-1} + \beta_i'' f_i +  \gamma_i'' f_{i+1} + \xi_i'' f_{i+2}, \\
\end{aligned}
\end{equation}
with
\begin{equation}\label{PLF11}
\begin{aligned}
\alpha_i' &:= \frac{-2 + 6w_i - 3w_i^2}{6 (\Delta \phi_0)}, && \alpha_i'' := \frac{1 - w_i}{(\Delta \phi_0)^2}, &&
\beta_i' := \frac{-1 - 4w_i + 3w_i^2}{2(\Delta \phi_0)}, && \beta_i'' := \frac{-2 + 3w_i}{(\Delta \phi_0)^2}, \\
\gamma_i' &:= \frac{2 + 2 w_i - 3w_i^2}{2(\Delta \phi_0)}, && \gamma_i'' := \frac{1 - 3w_i}{(\Delta \phi_0)^2}, && \xi_i' := \frac{-1 + 3 w_i^2}{6(\Delta \phi_0)}, && \xi_i'' := \frac{w_i}{(\Delta \phi_0)^2}. 
\end{aligned}
\end{equation}
One can then verify that $f_{\text{int},N}$ is $C^2$ on $[0,\phi_m]$. Furthermore, 
since $f$ and $g$ are smooth except at some $\phi_c$ where it is only $C^2$, one can show that $\{f''_{\text{int},N}\}_N$ is a globally equi-Lipschitz function on $[0,\phi_m]$ and hence equicontinuous. In addition, by \eqref{PLF7} and Taylor's theorem, we observe that $f_{\text{int},N}, f_{\text{int},N}'$ and $f''_{\text{int},N}$ converges to $f, f',$ and $f''$ pointwise. Since the family of functions $\{f''_{\text{int},N}\}_N$ is equibounded and equicontinuous, by Arzel\'a-Ascoli theorem, there exists a uniformly convergent subsequence $f''_{\text{int},N_k}$ of $f''_N$ under the sup norm. Repeating a similar argument and passing along subsequences for $f$ and $f'$, we can conclude that there exists a subsequence such that $f_{\text{int},N_k}$ converges to $f$ under the $C^2$ norm. A similar argument holds for $g$. Furthermore, since $K$ is compact, there exists a constant $C_3$ depending on $K$ such that $|h| \leq C_3$. Henceforth, there exists a positive integer $N$ such that 
\begin{equation}\label{PLF12}
|f''(\phi_i) - f''_{\text{int},N}(\phi_i)|,|g''(\phi_i) - g''_{\text{int},N}(\phi_i)| \leq \frac{\varepsilon_1}{6C_3^3}.
\end{equation}
Since $f''$, $f''_{\text{int},N}$, $g''$, and $g''_{\text{int},N}$ are uniformly continuous on the compact set $[0,\phi_m]$, there exists a $\delta > 0$ such that for all $|\phi - \tilde{\phi}| < \delta$, we have $|f''(\phi) - f''(\tilde{\phi})|,|f_{\text{int},N}''(\phi) - f_{\text{int},N}''(\tilde{\phi})|, |g''(\phi) - g''(\tilde\phi)|,|g_{\text{int},N}''(\phi) - g_{\text{int},N}''(\tilde{\phi})| < \frac{\varepsilon_1}{6C_3^3}$. Pick $\phi_i = \ml i + \frac{1}{2}\mr \delta$ for $i = 1, \cdots, \lfloor \frac{\phi_m}{\delta} - \frac{1}{2} \rfloor$. By a $3\varepsilon$ argument for both $f$ and $g$, one can deduce that 
\begin{equation}\label{PLF13}
\|f - f_{\text{int},N}\|_{C^2([0,\phi_m])}, 
\|g - g_{\text{int},N}\|_{C^2([0,\phi_m])} \leq \frac{\varepsilon_1}{2C_3^3},
\end{equation}
which in turn implies that
\begin{equation}\label{PLF14}
\|(F,G)-(F_{\text{int},N},G_{\text{int},N})\|_{\ml C^2(K)\mr^2} < \varepsilon_1.
\end{equation}
The conclusion then follows by applying Theorem \ref{Approxflux_generic}. \hfill $\square$

\subsection{Approximating Flux Functions}
Even though interpolating the flux functions via a cubic spline would work, this still requires a sufficiently fine grid over the parameter space. Moreover, this method scales exponentially in time with the number of parameters. In particular, for a bidisperse mixture with differing densities as in \cite{Bidensity}, the corresponding flux functions $f$ and $g$ are parameterized by $\phi_0 \in [0,\phi_m]$ and $X_0 \in [0,1]$, with the latter representing the fraction of particles to be of the first type. Hence, this warrants a more efficient numerical approximation of the flux functions and is the main motivation for Theorem \ref{Approxflux_generic}.

In particular, due to the functional form of the flux functions $F$ and $G$ depending only on functions of a single variable $f$ and $g$ defined on a compact set $\phi_0 \in [0,\phi_m]$, we can first approximate $f$ and $g$ uniformly in $C^2([0,\phi_m])$ with piecewise $C^2$ polynomials $f_p$ and $g_p$ respectively. These then correspond to the approximated flux functions $F_p(h,h\phi_0) = h^3 f_p\ml \frac{h\phi_0}{h}\mr$ and $G_p(h,h\phi_0) = h^3 g_p \ml \frac{h \phi_0}{h}\mr$ for $(h,h\phi_0) \in U$. For generic left and right state and any appropriate choice of a compact set $K$, since there exists a constant $C_3$ depending on the compact set $K$ such $|h| \leq C_3$, we then have the following estimates:
\begin{equation}\label{PLF6}
\begin{aligned}
\|F - F_p\|_{C^2(K)} \leq C_3^3 \|f - f_p\|_{C^2([0,\phi_m])}, \;
\|G - G_p\|_{C^2(K)} \leq C_3^3 \|g - g_p\|_{C^2([0,\phi_m])}. \\
\end{aligned}
\end{equation}
For a given $\varepsilon_1 > 0$, if pick $f_p$ and $g_p$ such that
$$ \|f - f_p\|_{C^2([0,\phi_m])},  \|g - g_p\|_{C^2([0,\phi_m])} \leq \frac{\varepsilon_1}{2C_3^3},$$ 
then $$\|(F,G) - (F_p,G_p)\|_{\ml C^2(K) \mr^2} < \varepsilon_1.$$
This then implies that Theorem \ref{Approxflux_generic} holds, implying generic structural stability under suitable polynomial approximations to the functions $f$ and $g$.

\subsubsection{Piecewise Polynomial Approximation} In this subsection, we will demonstrate the effectiveness of piecewise polynomial approximations for $f$ and $g$ in preserving the structural stability of Riemann solutions. To match most of the parameters in \cite{Murisic2011}, we set $\alpha = 25^\circ$, and hence the parameters in \eqref{PLF3} and \eqref{PLF4} are given by $\mu_l = 0.971$, $B_2 = 1.80036$, and $\phi_c = 0.50297$, with $\phi_m = 0.61$. As mentioned above, it can be verified numerically that $f$ and $g$ are $C^2$ across $\phi_c$. The comparison is done with functions $f$ and $g$ obtained from a cubic spline interpolation described in Section \ref{PLF_intp} with $\Delta \phi_0 = 0.001$. For our choice of piecewise $C^2$ polynomial approximations, we have
\begin{equation}\label{PLF_Numerics_1}
f_p(\phi_0) = \left\{ \begin{aligned}
\sum_{j=1}^{10} \beta_{f,j}^S (\phi_c - \phi_0)^{j-1} & \text{ for } \phi_0 < \phi_c, \\
\sum_{j=1}^{10} \beta_{f,j}^R (\phi_0 - \phi_c)^{j-1} & \text{ for } \phi_0 > \phi_c, \\
\end{aligned} \right.
\end{equation}
and
\begin{equation}\label{PLF_Numerics_2}
g_p(\phi_0) = \left\{ \begin{aligned}
\sum_{j=1}^{10} \beta_{g,j}^S (\phi_c - \phi_0)^{j-1} & \text{ for } \phi_0 < \phi_c, \\
\sum_{j=1}^{10} \beta_{g,j}^R (\phi_0 - \phi_c)^{j-1} & \text{ for } \phi_0 > \phi_c, \\
\end{aligned} \right.
\end{equation}
The choice of splitting the polynomials at $\phi_0 = \phi_c$ is motivated by the phase transition and lack of regularity across $\phi_c$. In addition, we will impose the following conditions on the coefficients of the polynomials as follows:
\begin{itemize}[leftmargin=30pt]
\item[(I)] $f_p$ and $g_p$ are $C^2$ across $\phi_0 = \phi_c$,
\item[(II)] $f_p(\phi_c) = \frac{\mu_l}{3}(1+ \rho_s \phi_c)\ml 1 - \frac{\phi_c}{\phi_m} \mr^2$, $g_p(\phi_c) = \phi_c f(\phi_c)$
\item[(III)] $f_p(0)  = \frac{\mu_l}{3}$, $g_p(0) = 0$,
\item[(IV)] $f_p(\phi_m) = g_p(\phi_m) = f_p'(\phi_m) = g_p'(\phi_m) = 0$.
\end{itemize}
The last two constraints were derived in \cite{Liwangshock}, and the second constraint at $\phi_c$ can be derived asymptotically using a similar argument as in \cite{Liwangshock}. Furthermore, we sample points $\phi_j$ on $[0,\phi_m]$ for which $f(\phi_j)$ and $g(\phi_j)$ are evaluated, according to the following rules:
\begin{itemize}[leftmargin=*]
\item 10 uniformly-spaced points are sampled close to $\phi_c$ corresponding to a phase transition, 
\item 10 uniformly-spaced points are sampled close to $\phi_m$ which corresponds to the suspension coalescing into a particle-rich ridge and exhibiting a different physical behavior, and
\item $13$ sets of 3 sparse points sampled across $[0,\phi_m]$. 
\end{itemize}
Note that points are sampled in sets of 3 so that centered finite differences can be used to obtain information about the derivatives of $f$ and $g$. Similarly, information about the derivatives for points sampled in the first two scenarios is provided for the interior points.

Recall that as observed in \cite{Liwangshock}, the solutions to \eqref{PLF1} exhibit singular shock for some left and right state, and this happens as $\phi_0$ approaches $\phi_m$ and $h$ approaches infinity. This thus justifies sampling extra points around $\phi_m$ and imposing constraints on $f$ and $g$ at $\phi_m$ to aid the polynomials in mimicking this behavior.

To obtain the coefficients $\beta$, we do so by minimizing the following loss function:
\begin{equation}\label{PLF_numerics_3}
\boldsymbol{\beta}_f = \frac{1}{2}\text{argmin}_{\boldsymbol{\beta}_f \in \mathbb{R}^{20}} \|\mathbf{f} - \mathbf{X} \boldsymbol{\beta}_f \|_2^2 + \frac{\lambda}{2} \|\mathbf{f}' - \mathbf{X}' \boldsymbol{\beta}_f \|_2^2 
\end{equation}
with $\lambda \geq 0$,
\begin{equation}\label{PLF_numerics_4}
\begin{aligned}
& \boldsymbol{\beta}_{f} = \begin{pmatrix} \boldsymbol{\beta}^S_f \\ \boldsymbol{\beta}^R_f
\end{pmatrix}, \text{ with } (\boldsymbol{\beta}^S_f)_j = \beta^{S}_{f,j}, (\boldsymbol{\beta}^R_f)_j =  \beta^{R}_{f,j} \quad \text{ for } j \in [[1,10]], \\
&(\mathbf{f})_j = f(\phi_j) \; \text{ for } j \in [[1,N]], \text{ and } (\mathbf{f}')_j = f'_{\text{int}}(\phi_j), \; \text{ for } j \in [[1,N']],\\
& \mathbf{X} = \begin{pmatrix}
\mathbf{X}^S & \mathbf{0}_{N_1 \times 10} \\ 
\mathbf{0}_{N_2 \times 10}  & \mathbf{X}^R \\ 
\end{pmatrix}, \, \mathbf{X}' = \begin{pmatrix}
\mathbf{X}'^S & \mathbf{0}_{N'_1 \times 10} \\ 
\mathbf{0}_{N'_2 \times 10}  & \mathbf{X}'^R \\ 
\end{pmatrix}, \, 
\end{aligned} 
\end{equation}
and 
\begin{equation}\label{PLF_numerics_5}
\begin{aligned}
&(\mathbf{X}^S)_{i,j} = (\phi_c - \phi_i)^{j-1}  
&& \text{ for } (i,j) \in [[1,N_1]] \times [[1,10]], \\
&(\mathbf{X}^R)_{i,j} = (\phi_i - \phi_c)^{j-1}  
&& \text{ for } (i,j) \in [[1,N_2]] \times [[1,10]], \\
&(\mathbf{X}'^S)_{i,j} =  -(j-1)(\phi_c - \phi_i)^{j-2}  &&\text{ for } (i,j) \in [[1,N_1']] \times [[1,10]], \\
&(\mathbf{X}'^R)_{i,j} = (j-1)(\phi_i - \phi_c)^{j-2}  &&\text{ for } (i,j) \in [[1,N_2']] \times [[1,10]], \\
\end{aligned}
\end{equation}
subject to (I), (II), (III), and (IV).
In \eqref{PLF_numerics_4} and \eqref{PLF_numerics_5} above, $N_1$ and $N_2$ represent the number of sampled points to the left and right of $\phi_c$ providing information on $f(\phi_i)$, $N_1'$ and $N_2'$ for the number of points with information on $f'(\phi_0)$ and $f''(\phi_0)$, with $N := N_1 + N_2$ and $N' := N_1' + N_2'$. The loss function above penalizes the distance between the values of $f_p$ and $f$ and their first derivative with $\lambda$ as the weight parameter. Note that we do not fit for the second derivative data as these data obtained via finite difference are not numerically stable, especially for a stiff nonlinear ODE as described in \eqref{PLF4}.
An analogous loss function for $g$ will be minimized too to obtain $\boldsymbol{\beta}_g$. Obtaining the optimal $\boldsymbol{\beta}_f$ from \eqref{PLF_numerics_3} will be quick with by utilizing appropriate interior point methods as this is a quadratic program with linear constraints.

\subsubsection{Numerical Results} To generate $f(\phi_i)$ and $g(\phi_i)$ at each sampled grid point $\phi_i$, we solve the nonlinear system of ODE in \eqref{PLF4} using an implicit Runge-Kutta method of the Radau IIA family of order 5 available from the scipy package in Python \cite{Radau}. The quadratic program in \eqref{PLF_numerics_3} is to be solved using the cvxopt package in Python \cite{cvxopt}. The script used for generating the figures below is available on \href{https://github.com/HK-Tan/Generic-Structural-Stability---Numerical-Simulations-for-Particle-Laden-Flow}{Github}.\footnote{Link: \href{https://github.com/HK-Tan/Generic-Structural-Stability---Numerical-Simulations-for-Particle-Laden-Flow}{https://github.com/HK-Tan/Generic-Structural-Stability---Numerical-Simulations-for-Particle-Laden-Flow}} We do so by fitting polynomials without derivative data ($\lambda = 0$) and with derivative data ($\lambda = 0.03$), with $\lambda$ obtained via a leave-one-out cross validation procedure.
\begin{figure}[htbp]
\centering
\includegraphics[width=0.47\linewidth]{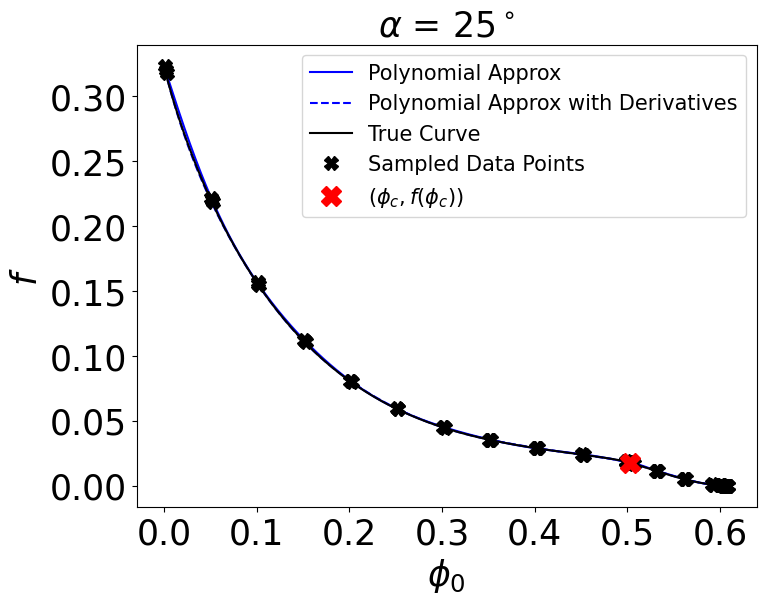}
\includegraphics[width=0.47\linewidth]{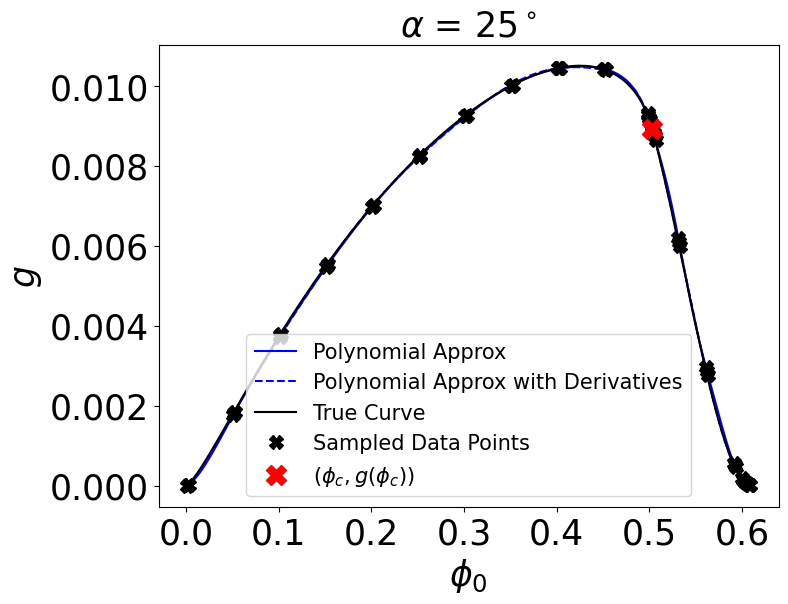} \\
\includegraphics[width=0.47\linewidth]{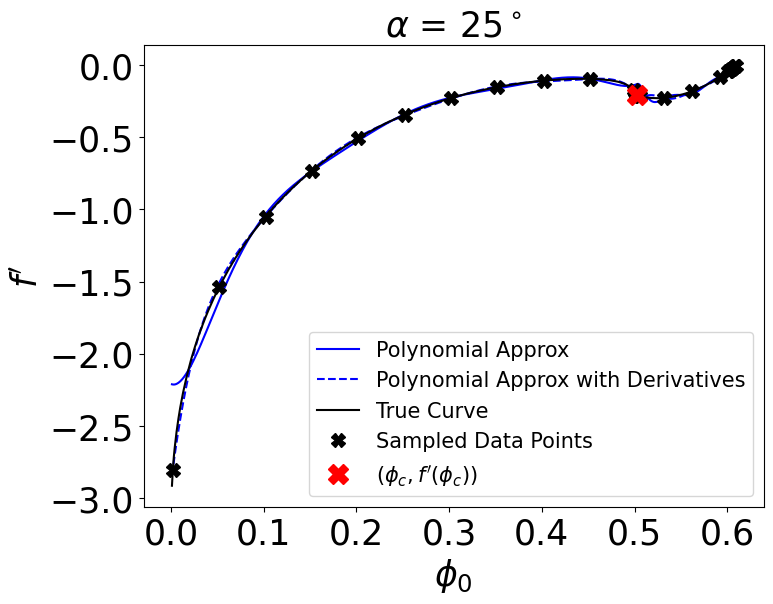}
\includegraphics[width=0.47\linewidth]{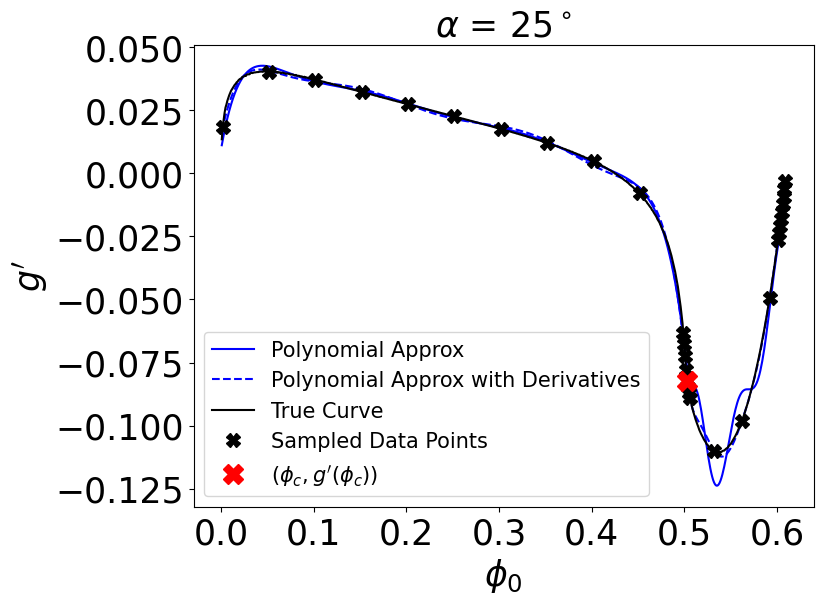} \\
\caption{Comparing the different types of polynomial approximation with the interpolation approximation for $f, f', g,$ and $g'$. The interpolated functions are labeled as `true curve' for each of the plots above.}
\label{fig:f,g,fp,gp}
\end{figure}
Figure \ref{fig:f,g,fp,gp} compares the fit between these and the interpolated functions $f$ and $g$. We can see from the plots that the piecewise polynomial approximations with or without derivative data approximate the functions well. However, for $f'$ and $g'$, the approximations without derivative data are not as good. This can be observed from the divergence close to $\phi_0 = 0$ for $f'$ and in the ridged regime for $g'$. Even though the fit for $g'$ could be improved with more points sampled in the ridged regime, the plots indicate that fitting with derivative data is still better for the same number of sparsely sampled points.

\begin{figure}[htbp]
\centering
\includegraphics[width=0.47\linewidth]{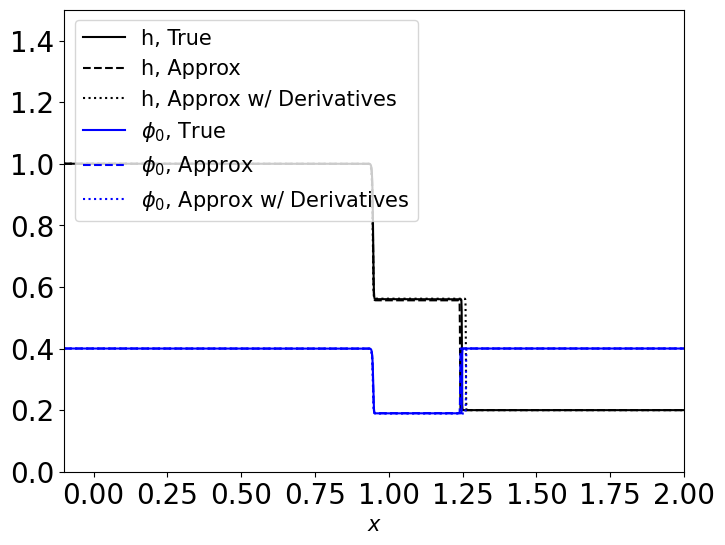} 
\includegraphics[width=0.47\linewidth]{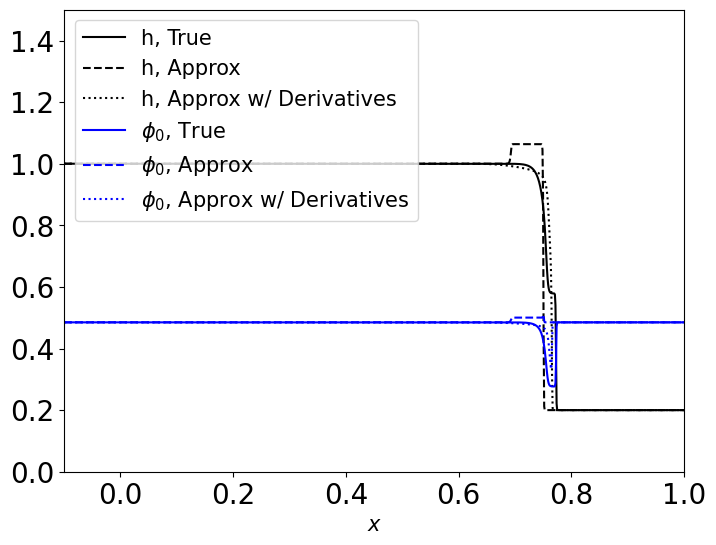} \\
\caption{Comparing the solution to the Riemann problem obtained via various flux approximations. Left: $(h_l,(\phi_0)_l) = (1,0.4)$ and $(h_r,(\phi_0)_r) = (0.2,0.4)$. Right: $(h_l,(\phi_0)_l) = (1,0.485)$ and $(h_r,(\phi_0)_r) = (0.2,0.485)$.}
\label{PLF:PDE}
\end{figure}

Next, we solve the system of conservation laws numerically for
\begin{itemize}[leftmargin=30pt]
\item[(1)] $(h_l,(\phi_0)_l) = (1.0,0.4)$ and $(h_r,(\phi_0)_r) = (0.2,0.4)$,
\item[(2)] $(h_l,(\phi_0)_l) = (1.0,0.485)$ and $(h_r,(\phi_0)_r) = (0.2,0.485)$.
\end{itemize}
In both scenarios, we use an upwind scheme with grid sizes $\Delta t = 0.0005$ and $\Delta x = 0.001$ up to time $t = 30$. The solutions obtained with the different types of approximation can be observed in Figure \ref{PLF:PDE}.
Note that the solutions obtained by the interpolated flux functions are treated as the true solution. 

In Figure \ref{PLF:PDE}, one can see that the double shock structure is preserved for case (1). However, this is not the case for the rarefaction-shock in case (2), as it degenerates into a double shock for the polynomial approximation fitted without derivative data. This is not surprising given that the quality of approximation is low for $\phi_0$ close to the ridged regime, as observed in the plot for $g'$ in Figure \ref{fig:f,g,fp,gp}. Furthermore, genuine non-linearity seems to fail numerically along $\phi_0$ close to $\phi_c$, which implies that it is possible for the Riemann problem to admit solutions with more than two waves.

\subsection{Comparing Physical Models} As mentioned in the introduction, the two competing classes of models used for gravity-driven particle-laden thin films are the diffusive flux and the suspension balance models. In \cite{ugrad_comparative} and \cite{Fastequilibrium}, the authors compared how the functions $f$ and $g$ obtained for two different classes of models affect the fluid velocity, the fluid front, and the height profile of the suspension. Nonetheless, the equilibrium profiles obtained hint at the functions $f$ and $g$ for both classes of models to be close. Let the superscript `DF' denote the functions $f$ and $g$ obtained under the diffusive flux models, and `SB' denote that obtained under the suspension balance models. For any given left and right state and the corresponding compact set $K$, there exists a constant $C_3$ depending on $K$ such that $|h| \leq C_3$. This implies that
\begin{equation}
\|(F^{\text{DF}},G^{\text{DF}}) - (F^{\text{SB}},G^{\text{SB}})\|_{\ml C^2(K) \mr^2} \leq C_3^3(\|(f^{\text{DF}},g^{\text{DF}}) - (f^{\text{SB}},g^{\text{SB}})\|_{\ml C^2(K) \mr^2}.
\end{equation}
Hence, if $\|(f^{\text{DF}},g^{\text{DF}}) - (f^{\text{SB}},g^{\text{SB}})\|_{\ml C^2(K) \mr^2} < \frac{\varepsilon_1}{C_3^3}$, we can apply Theorem \ref{Approxflux_generic} to conclude that the Riemann problem exhibits double-wave entropy solution of the same type and hence similar physical behavior for both models, with their intermediate states sufficiently close. 

\subsection{Measurement Errors}

\addt{In addition to the above applications, one can also use this theorem to argue that a similar physical structure as solutions to the Riemann problem \eqref{PLF1} and \eqref{PLFID} persists even in the presence of slight measurement errors to the left and right states. In particular, these measurement errors manifest as perturbations to the left and right states, in which we know from the application of Theorem \ref{Approxflux_generic} that for almost every left and right state, structural stability persists under such perturbations. This is useful especially given the fact that for conducting experiments similar to the Riemann problem as described above, it is usually difficult to measure the particle volume fraction $\phi_0$ of the corresponding slurry mixture.}

\section{Conclusions and Future Work}
\label{sec:future_work}

In this paper, we show that for generic left and right states, the Hugoniot loci and rarefaction curves intersect transversely at the intermediate state. Furthermore, we have shown that this is the exact requirement for the system to exhibit structural stability of Riemann solutions to the $2 \times 2$ system of hyperbolic conservation laws with \addt{respect to perturbations to the flux functions, the left state and the right state}. The argument generalizes easily to include structural stability with respect to perturbations to the left and right states.

Our above analysis focuses only on classical shocks singled out by the Lax entropy conditions \cite{Lax} and rarefaction connections.
Furthermore, we do not address more complex shock and rarefaction structures, such as composite waves, undercompressive (transitional) shocks, and contact discontinuities as described in \cite{Schecter}. These structures are more suitable to be analyzed from a vanishing viscosity perspective while utilizing a dynamical system argument to search for traveling wave connections between given left and right states. 

We also remark that it is not clear if the argument in this paper generalizes easily to $2 \times 2$ systems that do not satisfy the genuine non-linearity condition. One might have to investigate how we can incorporate Liu's entropy condition \cite{Liu2x2}, an analog to the Olenik chord condition for non-convex scalar conservation laws \cite{fronttracking}, into the arguments above. This is crucial as the system of conservation law system governing the evolution of the fluid mixture in \eqref{PLF} is not necessarily genuinely non-linear.

Even if we address all the points mentioned earlier, it was investigated in \cite{SurfaceTension2018} that a vanishing viscosity might not be the correct physical approach for the system in \eqref{PLF}. Instead, surface tension plays a crucial role, introducing a fourth-order regularizing term to the system of conservation laws in \eqref{PLF} rather than the usual second-order diffusion term associated with vanishing viscosity.

Last but not least, we note that it might be possible to generalize the arguments above to $3 \times 3$ or even $n \times n$ systems of conservation laws, at least in the class of $n$-wave solutions consisting of classical shocks, rarefactions, and possibly contact discontinuities singled out by Lax's entropy conditions. A potential transversality argument would involve the generic transverse intersection between a wave curve from a given state, and the ``foliated wave surface" emanating from the other state. If this is to be done using a coordinate-free approach,  we might be able to circumvent the graph condition as in (iii). This has applications to the compressible Euler equation and the generalization of \eqref{PLF} to a $3 \times 3$ system governing the evolution of gravity-driven particle-laden thin films with two different types of particles, such as particles with differing densities and sizes \cite{Bidensity, luong2025brazilnuteffectbidisperseparticle}.

\appendix
\section{Proof of Lemma \ref{C1maps}}\label{Proof_of_C1} 

To do so, we compute the relevant Fr\'echet derivatives as follows. Let $R > 0$ be positive real numbers such that the compact set $K$ lies in some open ball centered at $(0,0)$ of radius $R$.

Starting from $H: K \times \mathcal{O} \times U \rightarrow \mathbb{R}$ and its definition in \eqref{RH}, for any given $(F,G,u_g,v_g) \in \mathcal{O} \times U$, we have
\begin{equation}\label{C1lemma_1}
\begin{aligned}
& D_{(u,v)}H(\cdot;F,G,u_g,v_g) \\
=&\; \begin{pmatrix}
F_u(u,v) \cdot (v - v_g) - G_u(u,v) \cdot (u - u_g) - \ml G(u,v) - G(u_g,v_g)\mr \\
F_v(u,v) \cdot (v - v_g) - G_v(u,v) \cdot (u - u_g) + \ml F(u,v) - F(u_g,v_g)\mr \\
\end{pmatrix}^\intercal.
\end{aligned}
\end{equation}
Let $C$ be a positive real number such that 
$\max\{ \|F\|_{C^2(K)}, \|G\|_{C^2(K)} \} \leq C.$
From \eqref{C1lemma_1}, it is not difficult to see that for any $(u,v) \in K$,
\begin{equation}\label{C1lemma_2}
\begin{aligned}
\|D_{(u,v)}H(u,v;F,G,u_g,v_g)\|_{l^2} \leq 6C\|(u,v) - (u_g,v_g)\|_{l^2}.
\end{aligned}
\end{equation}
From \eqref{C1lemma_2}, we can deduce that
$D_{(u,v)}H(\cdot;F,G,u_g,v_g): K \rightarrow \mathcal{B}(K,\mathbb{R})$ is a bounded linear map. Next, for any given $(u,v,u_g,v_g) \in K \times U$, all $(F,G) \in \mathcal{O}$ and $(F_\delta,G_\delta) \in \ml C^2(K)\mr^2$, 
\begin{equation}\label{C1lemma_3}
\begin{aligned}
& D_{(F,G)}H(u,v;F,G,u_g,v_g)(F_\delta,G_\delta) \\
=&\; (F_\delta(u,v) - F_\delta(u_g,v_g))\cdot(v-v_g) - (G_\delta(u,v) - G_\delta(u_g,v_g))\cdot (u-u_g).
\end{aligned}
\end{equation}
Hence, we have the following estimate:
\begin{equation}\label{C1lemma_4}
\begin{aligned}
& |D_{(F,G)}H(u,v,F,G;u_g,v_g)(F_\delta,G_\delta)| \leq 4R\|(F_\delta,G_\delta)\|_{\ml C^2(K)\mr^2}.
\end{aligned}
\end{equation}
The above estimate implies that $D_{(F,G)}H(u,v;\cdot,u_g,v_g): \mathcal{O} \rightarrow \mathcal{B}(\ml C^2(K)\mr^2,\mathbb{R})$ is a bounded linear map. \addt{Last but not least, for any given $(u,v,F,G) \in K \times (C^2(K))^2$, we can use the symmetric property of the Hugoniot objective function as described in \eqref{T17} to obtain the exact same bounds as in \eqref{C1lemma_2} for $D_{(u_g,v_g)}H(u,v;F,G,u_g,v_g)$ and conclude that it is a bounded linear map from $\mathbb{R}^2$ to $\mathbb{R}$. Henceforth, $H: K \times \mathcal{O} \times U \rightarrow \mathbb{R}$ is a $C^1$ map.}

We now move on to the maps $R_k: K \times \mathcal{O} \times U \rightarrow \mathbb{R}$ for $k \in \{1,2\}$. For any given $(F,G,u_g,v_g) \in \mathcal{O} \times U$ and each $(u,v) \in K$, by \eqref{rareODE}, \eqref{rareODE2}, and \eqref{rareobj}, we have
\begin{equation}\label{C1lemma_5}
D_{(u,v)}R_k(u,v;F,G,u_g,v_g) = \begin{pmatrix}
- \Xi_k(u,v_k(u;F,G,u_g,v_g),F,G) & 1
\end{pmatrix}.
\end{equation}
Since $(u,v_k(u;F,G,u_g,v_g)) \in K$ and by the uni-directional assumption for $(F,G) \in \mathcal{O}$, we have that $D_{(u,v)}R_k(\cdot;F,G,u_g,v_g): K \rightarrow \mathcal{B}(K,\mathbb{R})$ is a bounded linear map. \addt{Similarly, for any given $(u,v,F,G) \in K \times (C^2(K))^2 $, we can utilize the derivative of the rarefaction objective functions with respect to the given state in \eqref{T22} and \eqref{T30} and the uni-directional assumption to conclude that $D_{(u,v)}R_k(u,v;F,G,\cdot): U\rightarrow \mathcal{B}(\mathbb{R}^2,\mathbb{R})$.}

\sloppy
To conclude this proof, we would like to show that the map $D_{(F,G)}R_k(u,v;\cdot,u_g,v_g): \mathcal{O} \rightarrow \mathcal{B}\ml \ml C^2(K)\mr^2, \mathbb{R} \mr$ is bounded. To do so, we first view $v_k(u;F,G,u_g,v_g)$ as a map $v_k: (u,F,G) \mapsto v_k(u;F,G,u_g,v_g)$ for a given state $(u_g,v_g)$, and hence $v_k(u,\cdot;u_g,v_g): \mathcal{O} \rightarrow \mathbb{R}$.
By the definition of $R_k$ in \eqref{rareobj}, this is equivalent to showing that the map
$D_{(F,G)}v_k(u;\cdot,u_g,v_g):\mathcal{O} \rightarrow \mathcal{B}\ml \ml C^2(K)\mr^2, \mathbb{R} \mr$ is bounded, with 
$$D_{(F,G)}R_k(u,v;\cdot,u_g,v_g) = D_{(F,G)}v_k(u;\cdot,u_g,v_g).$$
To compute this, we start from \eqref{T19},
dropping the dependence on $(u_g,v_g)$ while including the dependence of $\Xi_k$ on $(F,G)$ explicitly in the remaining parts of this proof for brevity. Next, we recall the chain rule for Fr\'echet derivatives as follows:
\begin{lemma}\label{chain} (Chain Rule.) Suppose that 
$f:U \subset X \rightarrow Y, \quad g: V \subset Y \rightarrow Z,$
where $U$ and $V$ are open subsets of $X$ and $Y$, respectively. If $f$ is differentiable at $x \in U$ and $g$ is differentiable at $f(x) \in V$, then $g \circ f$ is differentiable at $x$ with
$$d(g \circ f)(x) = dg(f(x)) \circ df(x).$$
\end{lemma}
Formally, one can apply Lemma \ref{chain} on \eqref{T19} such that for any $(F_\delta,G_\delta) \in C^2(K)^2$, we have
\begin{equation}\label{C1lemma_7}
\begin{aligned}
&\; D_{(F,G)}v_k(u,F,G)(F_\delta,G_\delta)  \\
= &\; \int^u_{u_g} D_v \Xi_k(s,v_k(s,F,G),F,G) D_{(F,G)}v_k(s,F,G)(F_\delta,G_\delta) \D s \\
&+ \int^u_{u_g} D_{(F,G)}\Xi_k(s,v_k(s,F,G),F,G)(F_\delta,G_\delta) \D s.
\end{aligned}
\end{equation}
Taking the derivative with respect to $u$ and applying the Fundamental Theorem of Calculus, we have
\begin{equation}\label{C1lemma_9}
\begin{aligned}
&\; \frac{\D}{\D u}D_{(F,G)}v_k(u,F,G)(F_\delta,G_\delta) \\
= &\; D_v \Xi_k(u,v_k(u,F,G),F,G) D_{(F,G)}v_k(u,F,G)(F_\delta,G_\delta)  \\
&+ D_{(F,G)}\Xi_k(u,v_k(u,F,G),F,G)(F_\delta,G_\delta).
\end{aligned}
\end{equation}
Let $Y(u) := D_{(F,G)}v_k(u,F,G)(F_\delta,G_\delta).$ Hence, \eqref{C1lemma_7} and \eqref{C1lemma_9} imply the following ODE:
\begin{equation}\label{C1lemma_10} \left\{
\begin{aligned}
Y'(u) &= A(u) Y(u) + B(u), \\
Y(u_g) &= 0,
\end{aligned}\right.
\end{equation}
with
\begin{equation}\label{C1lemma_11}
\begin{aligned}
A(u) &:= D_{v}\Xi_k(u,v_k(u,F,G),F,G), \\
B(u) &:= D_{(F,G)}\Xi_k(u,v_k(u,F,G),F,G)(F_\delta,G_\delta).
\end{aligned}
\end{equation}
The solution to the ODE in \eqref{C1lemma_10} is given by
\begin{equation}\label{C1lemma_12}
Y(u) = \int^u_{u_g} e^{\int^u_{\xi} A(s) \D s}B(\xi) \D \xi.
\end{equation}
To rigorously show that this is the correct Fr\'echet derivative, one can verify this using the definition or apply Dominated Convergence Theorem for interchanging Fr\'echet derivatives with an integral based on the estimates that we are able to compute.

Next, we compute an explicit expression for $A(v)$ and $B(v)$ above. For $A(v),$ by \eqref{rareODE2}, it is not difficult to see that 
\begin{equation}\label{C1lemma_13}
D_v \Xi_k(\cdot,F,G) := -\frac{F_v(\cdot)  (F_{uv}(\cdot) - D_v \lambda_k(\cdot,F,G)) - F_{vv}(\cdot) (F_u(\cdot) - \lambda_k(\cdot,F,G))}{F_v^2(\cdot)}.
\end{equation}
By the graph condition for $(F,G) \in \mathcal{O}$, this implies that $|F_v|$ is bounded from below by a positive constant on the compact set $K$. In addition, since
\begin{equation}\label{C1lemma_14}
\lambda_k(\cdot,F,G) := \frac{1}{2}\ml F_u(\cdot) + G_v(\cdot)
+ (2k-3) \sqrt{\Delta(\cdot,F,G)} \mr
\end{equation}
with $\Delta(\cdot,F,G)$ as the discriminant given by 
\begin{equation}\label{C1lemma_15}
\Delta(\cdot,F,G) = \ml F_u(\cdot) + G_v(\cdot)\mr^2 - 4\ml F_u(\cdot)G_v(\cdot) - F_v(\cdot)G_u(\cdot)\mr,
\end{equation}
together with the strict hyperbolicity of $(F,G)$ on $K$, one could also deduce that $D_v \lambda_k(\cdot,F,G)$ is bounded on $K$. Thus, there is a positive constant $C_1 > 0$ such that
$|D_v \Xi(u,v,F,G)| \leq C_1$. For valid $(s,v_k(s,F,G)) \in K$, these imply
\begin{equation}\label{C1lemma_19}
|D_v \Xi_k (s,v_k(s,F,G),F,G)| \leq C_1.
\end{equation}
On the other hand, one can compute the partial Fr\'echet derivative of $\Xi_k$ with respect to $(F,G)$, utilizing Lemma \ref{chain}, to obtain
\begin{equation}\label{C1lemma_16}
\begin{aligned}
& D_{(F,G)}\Xi_k(\cdot,F,G)(F_\delta,G_\delta) \\
=&\; -\frac{F_v(\cdot)\ml (F_\delta)_u(\cdot) - D_{(F,G)}\lambda_k(\cdot,F,G)(F_\delta,G_\delta)\mr - (F_\delta)_v(\cdot)\ml F_u - \lambda_k(\cdot,F,G)\mr}{F_v^2(\cdot)}
\end{aligned}
\end{equation}
for all $(F_\delta,G_\delta) \in \ml C^2(K) \mr^2$. Furthermore, we have
\begin{equation}\label{C1lemma_17}
\begin{aligned}
D_{(F,G)}\lambda_k(\cdot,F,G)(F_\delta,G_\delta) =&\; \frac{1}{2}\ml (F_\delta)_u(\cdot) + (G_\delta)_v(\cdot)\mr \\
&+ \frac{2k-3}{4}\frac{D_{(F,G)}\Delta(\cdot,F,G)(F_\delta,G_\delta)}{\sqrt{\Delta(\cdot,F,G)}}
\end{aligned}
\end{equation}
with
\begin{equation}\label{C1lemma_18}
\begin{aligned}
D_{(F,G)}\Delta(\cdot,F,G)(F_\delta,G_\delta) =&\; (2F_u + 2F_v - 4 G_v)(\cdot)(F_\delta)_u(\cdot) + 4 F_v(\cdot) (G_\delta)_u(\cdot) \\
&+  (2F_u + 2F_v + 4 G_u)(\cdot)(F_\delta)_v(\cdot)  + 4 F_u(\cdot) (G_\delta)_v(\cdot).
\end{aligned}
\end{equation}
By (i) on $K$, we deduce that $|\Delta(\cdot,F,G)|$ is bounded from below on the compact set $K$. Hence, there exists a positive constant $C_2$ such that 
$$|D_{(F,G)}\lambda_k(u,v,F,G)(F_\delta,G_\delta)| \leq C_2 \|(F_\delta,G_\delta)\|_{\ml C^2(K) \mr^2}.$$
Hence, from \eqref{C1lemma_16}, there exists a positive constant $C_3 > 0$ such that
\begin{equation}\label{C1lemma_20}
|D_{(F,G)}\Xi_k(\cdot,F,G)(F_\delta,G_\delta)| \leq C_3 \|(F_\delta,G_\delta)\|_{\ml C^2(K) \mr^2}.
\end{equation}
Plugging our estimates for $A(\cdot)$ and $B(\cdot)$ from \eqref{C1lemma_19} and \eqref{C1lemma_20} respectively into \eqref{C1lemma_12}, we have the following estimate:
\begin{equation}\label{C1lemma_21}
\begin{aligned}
|D_{(F,G)}v_k(u,F,G)(F_\delta,G_\delta)| \leq \frac{C_3}{C_1}\ml e^{C_1|v - v_g|} - 1\mr \|(F_\delta,G_\delta)\|_{\ml C^2(K)\mr^2}.
\end{aligned}
\end{equation}
We can then conclude that $D_{(F,G)}v_k(u,\cdot): \mathcal{O} \rightarrow \mathcal{B}(\ml C^2(K) \mr^2,\mathbb{R})$ is a bounded linear map. \addt{Thus, we have that $R_k: K \times \mathcal{O} \times U \rightarrow \mathbb{R}$ is a $C^1$ map for $k \in \{1,2\}$.}
\fussy

\section{Proof of Theorem \ref{foliated_parametric_transversality}}
\label{Proof_of_FPT}

To prove Theorem \ref{foliated_parametric_transversality}, we recall the Morse-Sard theorem which is used in the proof of Thom's parametric transversality theorem (Theorem \ref{parametric_transversality}):
\begin{theorem}\label{Sard}
(Morse-Sard Theorem. \cite{Hirsch}) Let $\mathcal{M}$ and $\mathcal{N}$ be $C^r$ manifolds (with $r \geq 1$) with dimensions $m$ and $n$ respectively and $f: \mathcal{M} \rightarrow \mathcal{N}$ be a $C^r$ map. If
$$r > \max\{0,m-n\},$$
then the set of critical values of $f$ is of measure zero in $\mathcal{N}$.
\end{theorem}

We then give the proof of Theorem \ref{foliated_parametric_transversality} below.\vspace{3pt}

\textit{Proof of Theorem \ref{foliated_parametric_transversality}.} Let $\mathcal{W} = F^{-1}(\mathcal{Z})$, and $\pi: \mathcal{XP} \rightarrow \mathcal{P}$ be the projection map given by
$$\pi(x,p) = p \quad \text{ for each }(x,p) \in \mathcal{XP}.$$ If $F_p \pitchfork \mathcal{Z}$ whenever $p \in \mathcal{P}$ is a regular value of the restriction of the projection map $\pi|_{\mathcal{W}}$, by Morse-Sard theorem, this happens for almost every $p \in \mathcal{P}$ and we are done. By the generalized regular value theorem in Theorem \ref{GPIT}, $F \pitchfork \mathcal{Z}$ from assumption 5 implies that the codimension of $\mathcal{W} = F^{-1}(\mathcal{Z})$ in $\mathcal{XP}$ is equals to the codimension of $\mathcal{Z}$ in $\mathcal{Y}$, and hence $\dim \mathcal{XP} - \dim \mathcal{W} = \dim \mathcal{Y} - \dim \mathcal{Z}$. In addition, $\pi|_{\mathcal{W}}$ is a $C^r$ map as $\mathcal{XP}$ and $\mathcal{P}$ are $C^r$ manifolds by assumption 2.  Thus, Morse-Sard theorem holds for $\pi|_\mathcal{W}: \mathcal{W} \rightarrow \mathcal{P}$ if 
$$\begin{aligned}
r &> \max\{0, \dim \mathcal{W} - \dim \mathcal{P}\}   = \max\{0, \dim \mathcal{XP} - \dim \mathcal{Y} + \dim \mathcal{Z} - \dim \mathcal{P}\} \\
&  = \max\{0, \dim \mathcal{X} + \dim \mathcal{P} - \dim \mathcal{Y} + \dim \mathcal{Z} - \dim \mathcal{P}\} \\
& = \max\{0, \dim \mathcal{Y} + \dim \mathcal{Z} - \dim \mathcal{X}\},
\end{aligned}$$
which is true by assumptions 1 and 2. It remains to show that $F_p \pitchfork \mathcal{Z}$ when $p \in \mathcal{P}$ is a regular value of $\pi|_{\mathcal{W}}$.

To do so, we first recall that $F \pitchfork \mathcal{Z}$ from assumption 5 implies that for each $(x,p) \in F^{-1}(\mathcal{Z}) \subset \mathcal{XP}$,
\begin{equation}\label{T7}
dF_{(x,p)}(T_{(x,p)}\mathcal{XP}) + T_{F(x,p)} \mathcal{Z} = T_{F(x,p)} \mathcal{Y}.
\end{equation}
By assumption 3, we have
\begin{equation}\label{T11}
dF_{(x,p)}(T_x \mathcal{X}_p \times T_p\mathcal{P}) + T_{F(x,p)} \mathcal{Z} \cong T_{F(x,p)} \mathcal{Y}.
\end{equation}
To prove that $F_p \pitchfork \mathcal{Z}$, we would want to show that for each $x \in F_p^{-1}(\mathcal{Z}) \subset \mathcal{X}_p$, we have
\begin{equation}\label{T8}
(dF_p)_{x}(T_x \mathcal{X}_p) + T_{F_p(x)} \mathcal{Z} = T_{F_p(x)} \mathcal{Y}.
\end{equation}
To do so, we first observe the following. Since $F_p: \mathcal{X}_p \rightarrow \mathcal{Y}$ is the restriction of the domain $F:\mathcal{XP} 
 \rightarrow \mathcal{Y}$ to $\mathcal{X}_p \times \{p\}$, for each $v_x \in T_x \mathcal{X}_p$ and by assumption 3, we have
\begin{equation}\label{T12}
dF_{(x,p)}(v_x,0) = (dF_p)_x (v_x).
\end{equation}

Starting from an arbitrary $x \in F_p^{-1}(\mathcal{Z})$, there exists a $z \in \mathcal{Z}$ such that $F_p(x) = z$, which in turn implies that $F(x,p) = z$ with $(x,p) \in \mathcal{W} = F^{-1}(\mathcal{Z})$.
To show that \eqref{T8} is true, assume that we are given a vector $v \in T_{F_p(x)}\mathcal{Y} = T_{z}\mathcal{Y}$. By \eqref{T7}, there exists vectors $v_{(x,p)} \in dF_{(x,p)}(T_{(x,p)}\mathcal{XP})$ and $v_z \in T_{F(x,p)}\mathcal{Z} = T_z\mathcal{Z}$ such that 
\begin{equation}\label{T13}
v_{(x,p)} + v_z = v.
\end{equation} 
Since $v_{(x,p)}\in dF_{(x,p)}(T_x \mathcal{X}_p \times T_p\mathcal{P})$, there exists a vector $(w,s) \in T_x \mathcal{X}_p \times T_p \mathcal{P}$ such that $v_{(x,p)} = dF_{(x,p)}(w,s).$ Subsequently, we consider the following two cases:
\begin{itemize}[leftmargin=*]
\item \ul{$s = 0 \in T_p\mathcal{P}$.} By \eqref{T12}, $v_{x,p} = dF_{(x,p)}(w,0) = (dF_p)_x(w)$, and we have that $v_{x,p} \in (dF_p)_x(T_x \mathcal{X}_p)$ and $v_z \in T_z \mathcal{Z}$ such that $v_{x,p} + v_z = v$ by \eqref{T13}.
\item \ul{$s \neq 0 \in T_p\mathcal{P}$.} Since $p$ is a regular value of the restriction of the projection map $\pi|_{\mathcal{W}}$ and $d\pi|_{\mathcal{W}}$ is just the projection to the second vector, then there exists $u \in T_x \mathcal{X}_p$ such that $(u,s) \in T_{(x,p)} \mathcal{W}$ and $(d\pi|_\mathcal{W})_{(x,p)}(u,s) = s$. Note that elements in $T_{(x,p)}\mathcal{W}$ are in $T_{(x,p)}(\mathcal{XP}) \cong T_x \mathcal{X}_p \times T_p \mathcal{P}$ since $\mathcal{W}$ is a submanifold of $\mathcal{XP}$. Henceforth, as $F|_\mathcal{W}: \mathcal{W} \rightarrow \mathcal{Z}$, then $(dF|_\mathcal{W})_{(x,p)}: T_{(x,p)}\mathcal{W} \rightarrow T_z \mathcal{Z}$. Since $(u,s) \in T_{(x,p)}\mathcal{W}$, then $dF_{(x,p)}(u,s) \in T_z \mathcal{Z}$ a. From \eqref{T12}, this then implies that
\begin{equation}\label{T14}
dF_{(x,p)}((w,s) - (u,s)) = dF_{(x,p)}(w-u,0) = (dF_p)_x(w-u).
\end{equation}
We claim that $(dF_p)_x(w-u) \in (dF_p)_x(T_x \mathcal{X}_p)$ and $dF_{(x,p)}(u,s) + v_z \in T_z\mathcal{Z}$ such that their vector sum is equals to $v$. Indeed, by \eqref{T13} and \eqref{T14}, observe that
\begin{equation}\label{T15}
\begin{aligned}
&\underbrace{(dF_p)_x(w-u)}_{\in (dF_p)_x(T_x \mathcal{X}_p)} +  \underbrace{dF_{(x,p)}(u,s) + v_z}_{\in T_z\mathcal{Z}} \\
=&\; dF_{(x,p)}(w,s) - dF_{(x,p)}(u,s) +  dF_{(x,p)}(u,s) + v_z = v_{(x,p)} + v_z = v.
\end{aligned}
\end{equation}
\end{itemize}
Hence, \eqref{T8} is true in both cases, and we are done.
$\hfill \square$

\section*{Acknowledgments}
We thank Stephen Schecter for helpful conversations about his previous work in this area.  The first author would like to thank Sarah Burnett for helpful discussions regarding the particle-laden flow experiment at UCLA. The first author would also like to thank Joshua Enwright for helpful discussions on differential topology and transversality.


\bibliographystyle{siamplain}
\bibliography{references}
\end{document}